\documentclass[12pt]{article} 
\usepackage[sectionbib]{natbib}
\usepackage{array,epsfig,fancyheadings,rotating}
\usepackage[]{hyperref}  
\usepackage{sectsty, secdot}
\usepackage{amsmath}
\usepackage{graphicx,psfrag,epsf}
\usepackage{enumerate}
\usepackage{natbib}
\usepackage{url}
\usepackage{amsmath,amssymb,chemarrow}
\usepackage{graphicx}
\usepackage{subfigure}
\usepackage{colortbl,color,dcolumn}
\usepackage{todonotes}
\usepackage{enumerate}
\usepackage{eqnarray,array}
\usepackage{booktabs}
\usepackage{listings}
\usepackage{bm}
\usepackage{amsthm}
\usepackage{amsfonts}
\usepackage{subfigure}
\usepackage{caption}
\usepackage{graphicx}
\usepackage{wrapfig}
\usepackage{lscape}
\usepackage{rotating}
\usepackage{epstopdf}
\usepackage{lineno}

\sectionfont{\fontsize{12}{14pt plus.8pt minus .6pt}\selectfont}
\renewcommand{\theequation}{\thesection\arabic{equation}}
\subsectionfont{\fontsize{12}{14pt plus.8pt minus .6pt}\selectfont}

\usepackage{amsmath}
\usepackage{amssymb}
\usepackage{amsfonts}
\usepackage{multirow}
\usepackage{amsthm}

\newtheorem{theorem}{Theorem}
\newtheorem{lemma}{Lemma}
\newtheorem{corollary}{Corollary}

\theoremstyle{definition}

\newtheorem{example}{Example}

\begin{document}


\renewcommand{\baselinestretch}{2}

\renewcommand{\thefootnote}{}
$\ $\par


\fontsize{12}{14pt plus.8pt minus .6pt}\selectfont \vspace{0.8pc}
\centerline{\large\bf DESIGN BASED INCOMPLETE U-STATISTICS}
\vspace{.4cm} \centerline{Xiangshun Kong$^1$, Wei Zheng$^2$} \vspace{.4cm} \centerline{\it
    $^1$Beijing Institute of Technology and $^2$University of Tennessee} \vspace{.55cm} \fontsize{9}{11.5pt plus.8pt minus
.6pt}\selectfont


\begin{quotation}
\noindent {\it Abstract:}
U-statistics are widely used in fields such as economics, machine learning, and statistics. However, while they enjoy desirable statistical properties, they have an obvious drawback in that the computation becomes impractical as the data size $n$ increases. Specifically, the number of combinations, say $m$, that a U-statistic of order $d$ has to evaluate is $O(n^d)$. Many efforts have been made to approximate the original U-statistic using a small subset of combinations since Blom (1976), who referred to such an approximation as an incomplete U-statistic. To the best of our knowledge, all existing methods require $m$ to grow at least faster than $n$, albeit more slowly than $n^d$, in order for the corresponding incomplete U-statistic to be asymptotically efficient in terms of the mean squared error. In this paper, we introduce a new type of incomplete U-statistic that can be asymptotically efficient, even when $m$ grows more slowly than $n$. In some cases, $m$ is only required to grow faster than $\sqrt{n}$. Our theoretical and empirical results both show significant improvements in the statistical efficiency of the new incomplete U-statistic. 

\vspace{9pt}
\par
\end{quotation}\par

\def\thefigure{\arabic{figure}}
\def\thetable{\arabic{table}}

\renewcommand{\theequation}{\thesection.\arabic{equation}}

\fontsize{12}{14pt plus.8pt minus .6pt}\selectfont
\newpage
\setcounter{section}{1} 
\setcounter{equation}{0} 

\lhead[\footnotesize\thepage\fancyplain{}\leftmark]{}\rhead[]{\fancyplain{}\rightmark\footnotesize\thepage}
\noindent{\large\bf 1. Introduction}
\label{sec:intro}

The U-statistic has been a fundamental statistical estimator since the work of \cite{hoeffding:1948}, who studied its theoretical properties and established central limit theorems for non-degenerate U-statistics. \cite{eagleson:1979} derived asymptotic distributions of some degenerate U-statistics of order two, which were then extended to all degenerate U-statistics by \cite{lee:1979}. Other extensions include a variant of U-statistics called V-statistics by \cite{richard:1947}, U-statistics for stationary processes by \cite{enqvist:1985}, and multi-sample U-statistics by \cite{lehmann:1951} and \cite{sen:1974,sen:1977}.

The theory of U-statistics admits a minimum variance unbiased estimator of an estimable parameter for a large class of probability distributions, hence its popularity in applications. However, U-statistics can also be time consuming to compute. For a U-statistic of order $d$, the number of combinations, say $m$, to be evaluated is $\binom{n}{d}$, that is $O(n^d)$, where $n$ is the data size. Suppose $n=10^4$ and $d=3$. Then, listing the $\binom{10^4}{3}$ combinations requires $667$ GB of memory and a computing time of approximately $100$ hours on a Macbook Pro with Intel Core i7 2.9 GHz CPU. With $n=10^5$ and $d=4$, the required memory is roughly  $16.7$ EB and the computing time is projected to be $285,000$ years. To provide context, \cite{hilbert:2011} estimated that humankind was able to store 295 EB of optimally compressed data in 2007. 
The issue of computational difficulty becomes even more severe in the bootstrap approximation of the asymptotic distribution of a U-statistic; see, for instance, \cite{bickel:1981}, \cite{bretagnolle:1983}, \cite{dehling:1994}, and \cite{marie:1993a, marie:1993b}, among others.

For certain U-statistics, the computational complexity can be reduced to $O(n)$ by exploiting the structure of the kernel function, especially when the data are univariate and consist of one sample. However, in practice, such a computational reduction is often not feasible. Note that we do not focus here on which U-statistics are candidates for a reduction in the original computational complexity of $O(n^d)$ because our goal is to study a generic scheme for the fast approximation of U-statistics. A natural remedy is to take a sample of size $m\ll\binom{n}{d}$ from all possible combinations. \cite{blom:1976} referred to the resulting estimator as an {\it incomplete U-statistic}. The problem of identifying a good incomplete U-statistic is related to the design of the sampling scheme. Of the various options, the vanilla scheme of simple random sampling by \cite{blom:1976} has received much attention in the literature. \cite{janson:1984} established the asymptotic distributions of incomplete U-statistics based on random sampling (ICUR), \cite{herrndorf:1986} established the invariance principle for the statsitics, and \cite{chen:2018} studied the vector- and matrix-valued ICUR. For a more detailed discussion on incomplete U-statistics, refer to \cite{wang:2012} and \cite{wang:2014}.

First, we introduce some required notation. For $\alpha>0$, we use $m\prec n^{\alpha}$, $m\asymp n^{\alpha}$, and $m\succ n^{\alpha}$ to mean $m/n^{\alpha}\rightarrow 0$, $0<\varliminf m/n^{\alpha}\leq\varlimsup m/n^{\alpha}<\infty$, and $m/n^{\alpha}\rightarrow \infty$, respectively. For a given incomplete U-statistic, say $U$, its efficiency is defined in terms of the mean squared error (MSE): ${\rm Eff}(U)={\rm MSE}(U_0)/{\rm MSE}(U)$, where $U_0$ is the complete U-statistic. An incomplete U-statistic is said to be {\it asymptotically efficient} if ${\rm Eff}(U)\rightarrow 1$ as $n\rightarrow\infty$. Note that the ICUR is asymptotically efficient for the non-degenerate case when $m\succ n$; see (\ref{eq:eff}) for a theoretical verification, and Table $1$ for empirical evidence. 

Blom (1976) also proposed sampling schemes based on the design of an experiment. In particular, balanced incomplete block designs (BIBDs) have been examined by \cite{brown:1978} and \cite{lee:1982}. The latter also proved that incomplete U-statistics based on BIBDs achieve the minimum variance among all unbiased estimators for a given $m$. By Raghavarao (1971), a BIBD exists whenever $n=6a+3$ for any positive integer $a$. Unfortunately, the optimality of the BIBD does not make it practically attractive because its construction requires $m\asymp n^2$; see Table \ref{tb:introduction}. The same issue exists for the permanent design of \cite{rempala:2003} and the rectangular design of \cite{rempala:2004}. For the case of $m/n\rightarrow 1$, Blom (1976) proposed using a Latin square and a Graeco-Latin square to guide the sampling scheme. However, the efficiency of the estimator derived in this way is essentially asymptotically the same as that of the ICUR. Moreover, the limit of the efficiency does not exceed $d/(1+d)$ as $n\rightarrow \infty$; see (\ref{eq:eff}) and the follow-up discussion. 

Another method recently proposed in the literature is the divide and conquer (DC) strategy of \cite{lin:2010}, which randomly divides the data into many groups, calculates the complete U-statistic within each group, and then takes the average of these complete U-statistics. Unfortunately, the DC is even less efficient than the ICUR. Moreover, it is not available when $m\leq n$; see Table 1.

We conclude that the ICUR is still the most viable of the existing choices of incomplete U-statistics. It performs as well as a design-based method when a design exists. It also possesses several advantages, such as a flexible choice of $m$, the availability of asymptotic properties, and being extendable to multi-sample cases. 

In this paper, we introduce a new type of incomplete U-statistic that is substantially more efficient than the ICUR, while maintaining the latter's aforementioned advantages. It has three main steps: ($i$) Divide the data into $L (\ll n)$ groups of homogeneous units. $(ii)$ Judiciously select a collection of the combinations of the groups based on a design structure called an orthogonal array (OA). $(iii)$ Randomly select a combination of inputs from each selected group combination. We call the derived estimator the incomplete U-statistic based on division and an orthogonal array (ICUDO). Our first example provides a snapshot of the performance of the major incomplete U-statistics mentioned so far. 

\begin{example}\label{example:manycomparison} (The symmetry of distribution).
The kernel function $g(x_1,x_2,x_3)={\rm sign}(2x_1-x_2-x_3)+{\rm sign}(2x_2-x_1-x_3)+{\rm sign}(2x_3-x_1-x_2)$ has mean zero when the distribution of the data is symmetric. The data consists of $n=10^3$ independent and identically distributed (i.i.d.) observations generated iid from the standard normal distribution. The performance of the ICUR, BIBD, DC, and ICUDO is measured by their efficiency at different values of $m$. 
\end{example}

\begin{table}[h]
  \begin{center}
    \caption{Comparison of efficiencies in Example \ref{example:manycomparison}.}\label{tb:introduction}
    \begin{tabular}{|c|c|c|cccc|}\hline    
 $m$&$m/n$&$m/\binom{n}{3}$&ICUR&BIBD&DC&ICUDO\\\hline
 $1.0\times10^3$&$1.0$&$6.018\times10^{-6}$&21.62\%&2.706\%&--&36.31\%\\
 $1.2\times10^4$&$12.0$&$7.222\times10^{-5}$&74.97\%&9.155\%&55.60\%&100\%\\
 $5.7\times10^4$&$57.0$&$3.430\times10^{-4}$&97.40\%&21.81\%&76.70\%&100\%\\
 $1.66\times 10^5$&$116.0$&$1.000\times10^{-3}$&100\%&100\%&84.22\%&100\%\\
 $3.92\times10^5$&$392.0$&$2.359\times10^{-3}$&100\%&100\%&90.71\%&100\%\\
 $1.617\times 10^6$&$1617.0$&$9.731\times10^{-3}$&100\%&100\%&95.64\%&100\%\\
   \hline
    \end{tabular}
  \end{center}

\end{table}
Note that the DC is unavailable when $m\leq n$, and the BIBD does not exist in most cases, except for $m=166167$. For $m\leq 166167$, the sample size is separately reduced for the BIBD in order to make it available. The ICUR has the same efficiency as the BIBD method at $100\%$ when the BIBD exists. It is more efficient than the DC method whenever the DC is available. However, the ICUDO methods outperforms the ICUR for all $m$. 


Here, we briefly explain why our ICUDO performs so well. Note that existing design-based methods focus on the arrangement of indices of units, without referring to their actual values. The ICUDO method exploits the fact that replacing a unit by another one with a similar value does not change the value of the kernel function $g$ too much. For example, suppose the first six numbers of the data are $(1,2,3,1,2,3)$. Then, a kernel function of order three yields the same value by evaluating the first three and the next three units. Beyond the grouping idea, we use the OA to achieve the projective uniformity of the group combinations in the dominating lower-dimensional spaces. This allows us to recover information on the lower dimension's variability in the U-statistics, which is the dominating part of Hoeffding's decomposition of the U-statistics. As shown later, in the non-degenerate case, whereas the ICUR needs $m\succ n$ to be  asymptotically efficient, the ICUDO requires a substantially smaller $m$; sometimes even $m\succ \sqrt{n}$ will suffice. See Theorem \ref{thm:lipschitzonesample} for the latter case. When the U-statistic is degenerate, both methods require larger $m$, but the ICUDO still requires a substantially smaller $m$ than that of the ICUR.

The rest of the paper is organized as follows. Section 2 introduces the construction of the ICUDO for univariate data and derives its asymptotic properties. Section 3 discusses the debiasing issues of the ICUDO for the degenerate case. Section 4 constructs a debiased ICUDO for multi-dimensional data. Simulations are presented in each section to support the theoretical results. Section 5 concludes the paper and points out some future research topics. All proofs are postponed to the Appendix. Additional theorems are given in the online Supplementary Material. 
\vspace{0.5cm}
\setcounter{section}{2} 
\setcounter{equation}{1} 

\lhead[\footnotesize\thepage\fancyplain{}\leftmark]{}\rhead[]{\fancyplain{}\rightmark\footnotesize\thepage}
\noindent{\large\bf 2. ICUDO based on univariate data}
\label{sec:intro}

Let $X_1,\ldots,X_n$ be a random sample of size $n$ from a univariate distribution, say $F$. For a given symmetric kernel function, say $g: R^d\rightarrow R$, of order $d$, the uniformly minimum variance unbiased estimator (UMVUE) of the parameter $\Theta=\int g(x_1,\ldots,x_{d})$ $dF(x_1)\ldots dF(x_{d})$
is given by the U-statistic
\begin{equation}\label{eq:ustat}
U_0=\binom{n}{d}^{-1}\sum_{\bm \eta\in S_{n,d}}g({\mathcal{X}}_{\bm \eta}),
\end{equation}
where $S_{n,d}=\{{\bm \eta}=(\eta_1,\ldots,\eta_{d}): 1\leq \eta_1<\eta_2<\ldots<\eta_{d}\leq n\}$ and ${\mathcal{X}}_{\bm \eta}=(X_{\eta_1},\ldots,X_{\eta_d})$. When $S_{n,d}$ is replaced with the set of all $n^d$ ordered combinations, the corresponding average in (\ref{eq:ustat}) is called a {\it V-Statistic} (\cite{richard:1947}). The main difference is that V-statistics include combinations with duplicated units, such as $(1,1,2)$. Throughout this paper, we adopt the mild assumption $Eg^2\left(X_{{1}},\ldots,X_{{d}}\right)<\infty.$

Unless there is some special structure of $g$ that can be exploited to reduce the computational burden, in general, (\ref{eq:ustat}) becomes impractical to compute as $n$ increases. To address this problem, Blom (1976) proposed using the following incomplete U-statistic as a fast approximation:
\begin{equation}\label{eq:incompleteuonesample}
U=\frac{1}{m}\sum_{{\bm\eta}\in S}g(\mathcal{X}_{\bm\eta}),
\end{equation}
where $S\subset S_{n,d}$, with its cardinality $m=|S|$ being only a fraction of $\binom{n}{d}$. The statistic in (\ref{eq:incompleteuonesample}) becomes an ICUR when $S$ is a simple random sample, which we denote as $U_{\rm RND}$. 

Here, we briefly review the properties of $U_0$ and $U_{\rm RND}$. For arbitrary positive integers $N$ and $p$, define $\mathcal{Z}_N=\{1,\ldots,N\}$ and $\mathcal{Z}_N^p=\{(z_1.\ldots,z_p):z_j\in\mathcal{Z}_N, 1\leq j\leq p\}$. Following \cite{hoeffding:1948}, for ${\bm u}\subseteq \mathcal{Z}_{d}$ and ${\bm x}=(x_1,\ldots,x_d)$, denote
$g_{\bm u}({\bm x})=\int g({\bm x})dF_{{\bm u}^c}$,
with ${\bm u}^c=\mathcal{Z}_{d}\setminus{\bm u}$ and $dF_{{\bm u}}=\prod_{j\in{\bm u}}dF(x_j)$. With the conventions $g_{\emptyset}({\bm x})=\Theta$ and $h_{\emptyset}({\bm x})=0$, we recursively define the projection
\begin{equation}\label{eq:hufunctiononesample}
h_{\bm u}({\bm x})=g_{\bm u}({\bm x})-\sum_{\bm v\subseteq\mathcal{Z}_d:{\bm v}\subset{\bm u}}h_{\bm v}(x).\notag
\end{equation}
Because $g$ is symmetric, we have $Eg_{\bm v}^2=Eg_{\bm u}^2$ and $Eh_{\bm v}^2=Eh_{\bm u}^2$ for any pair ${\bm u},{\bm v}\subseteq \mathcal{Z}_{d}$, with $|{\bm v}|=|{\bm u}|$. Hence, we can now define
\begin{equation}\label{eq:varelementsonesample}
\sigma_{{j}}^2={\rm Var}(g_{\bm u})~{\rm and}~\delta_{{j}}^2={\rm Var}(h_{\bm u}), ~ ~{\rm with} ~|{\bm u}|=j.\notag
\end{equation}
Following \cite{hoeffding:1948} and \cite{blom:1976}, we have
\begin{eqnarray}
{\rm MSE}(U_0)&=&\binom{n}{d}^{-1}\sum_{j=1}^{d}\binom{d}{j}\binom{n-d}{d-j}\sigma_{j}^2~=~\sum_{j=1}^d\binom{d}{j}^2\binom{n}{j}^{-1}\delta_j^2,\label{eq:6}\\
{\rm MSE}(U_{\rm RND})&=&{\rm MSE}(U_0)+\frac{\sigma_d^2}{m}+O\left(\frac{1}{nm}\right)\notag\\
&=&{\rm MSE}(U_0)+\frac{1}{m}\sum_{j=1}^d\binom{d}{j}\delta_{j}^2+O\left(\frac{1}{nm}\right).~~~~\label{eq:generalchangeonesample}
\end{eqnarray}
In (\ref{eq:6}) and (\ref{eq:generalchangeonesample}), the MSEs are expressed in terms of both $\sigma_j^2$ and $\delta_j^2$. The equivalences are established by $\sigma_{j}^2=\sum_{j'=1}^j\binom{j}{j'}\delta_{j'}^2$, for $1\leq j\leq d$. The U-statistic and the kernel function $g$ are called {\it non-degenerate} if $\delta_1^2=\sigma_1^2>0$, and are called {\it order-q degenerate} if $\sigma_q^2=0$ and $\sigma_{q+1}^2>0$, or equivalently $\delta^2_1=\cdots=\delta^2_q=0$ and $\delta^2_{q+1}>0$. For the non-degenerate case, we have ${\rm Var}(U_0)\asymp n^{-1}$, which together with (\ref{eq:generalchangeonesample}) yields
\begin{equation}\label{eq:eff}
{\rm Eff}(U_{\rm RND})=\left\{
\begin{array}{rcl}
1-O(n/m), & & {m\succ n}\\
\frac{1}{1+\frac{n}{m}\frac{\sigma_d^2}{d^2\delta_1^2}}+O(1/n), & & {m\asymp n}\\
O(m/n), & & {m\prec n}.
\end{array} \right.
\end{equation}
As a result, we have ${\rm Eff}(U_{\rm RND})\rightarrow 1$ when $m\succ n$, ${\rm Eff}(U_{\rm RND})\rightarrow 0$ when $m\prec n$, and ${\rm Eff}(U_{\rm RND})\rightarrow \left(1+\frac{\sigma_d^2}{cd^2\sigma_1^2}\right)^{-1}$ when $m/n\rightarrow c$, for a constant $c>0$. With $c=1$, \cite{blom:1976} proposed using Latin squares and Graeco-Latin squares to construct the incomplete U-statistics. In such a case, we can verify that its efficiency is asymptotically the same as that of $U_{\rm RND}$, and $\lim_{n\rightarrow\infty}{\rm Eff}(U_{\rm RND})\leq d/(1+d)$, from (\ref{eq:eff}) and $\sigma_d^2\leq d\sigma_1^2$. In contrast, Theorem \ref{thm:noassumptiononesample} shows that the ICUDO is asymptotically efficient when $m\asymp n$. Stronger results are stated in Theorem \ref{thm:lipschitzonesample} in Section 2.1 and in similar theorems in the Supplementary Material under various conditions on $g$ and $F$.

\noindent{\bf 2.1. One-sample U-statistics}

Recall that $\delta_j^2={\rm Var}(h_{\bm{u}})$, for $|\bm{u}|=j$, $1\leq j\leq d$, and note that the coefficient of $\delta_j^2$ in (\ref{eq:6}) is $O(n^{-j})$. Hence, it is more important to capture the variability of $g$ in its lower-dimensional projected space. This idea matches perfectly with the projective property of the OA. An OA denoted by $OA(m,d,L,t)$, is an $m$ by $d$ array with entries from $\{1,\ldots,L\}$, arranged in such a way that for any $m$ by $t$ subarray, all ordered $t$-tuples of the entries from $\{1,\ldots,L\}$ appear $\lambda=m/L^t$ times in the rows. The number $t$ is called the strength of the OA; see the matrix $A$ defined in (\ref{eqn:210}) as an example of $OA(9,4,3,2)$. In this case, the ordered $2$-arrays are $\{(i_1,i_2):1\leq i_1,i_2\leq 3\}$. Consider any two columns of $A$, we can see that all these ordered $2$-tuples appear once, that is, $\lambda=1$.
For sets $\mathcal{S}_1,\ldots,\mathcal{S}_q$, define 
$\prod_{i=1}^q\mathcal{S}_i=\{(\bm{s}_1,\ldots,\bm{s}_q):\bm{s}_i\in\mathcal{S}_i\}$. The ICUDO is constructed as follows. For ease of illustration, we assume $n$ is a multiple of $L$. Actually, throughout the manuscript, we assume that $L\ll n$. Thus, we may randomly draw an $n'=\lfloor n/L\rfloor \cdot L$ subsample as the new data set. The information loss in this process is negligible compared with the original size $n$.

\begin{itemize}
\item[Step 1.]  Let $A_0$ be an $OA(m,d,L,t)$. Apply random level permutations $\{\pi_1,\ldots,\pi_d\}$ to columns of $A_0$ independently. Specifically, for $l\in\mathcal{Z}_L$, change all elements $l$ in the $j$th column of $A_0$ to $\pi_j(l)$. The new OA is denoted by $A=(a_{ij})_{m\times d}$.
\item[Step 2.] Create the partition $\mathcal{Z}_{n}=\bigcup_{l=1}^LG_l$ such that $|G_l|=n/L$ for $l\in\mathcal{Z}_L$, and $X_{i_1}\leq X_{i_2}$ for any $i_1\in G_{l_1}, i_2\in G_{l_2}$, with $l_1<l_2$.
\item[Step 3.] For $i=1,\ldots,m$, independently draw an element, say ${\bm \eta}^i$, uniformly from $\prod_{j=1}^{d}G_{a_{ij}}$.
the ICUDO based on the OA $A$ is defined as
\begin{equation}\label{uoa}
U_{oa}=\frac{1}{m}\sum_{i=1}^mg(\mathcal{X}_{\bm\eta^i}).
\end{equation}
\end{itemize}

The level permutation in step 1 ensures that each row of $A$ takes each $d$-tuple with equal probability. At the same time, the projective uniformity of the beginning OA, $A_0$, carries over to $A$. Here, we ensure that $A$ is free of a coincidence defect, which means no two rows are the same in any $m\times(t+1)$ subarray. This property is necessary for the relevant theorems to hold. Step 2 divides the data into homogeneous groups. Step 3 is built on the first two steps. It chooses representative elements from selected groups, and the selection of groups is guided by the structure of $A$. Note (\ref{uoa}) is in the form of (\ref{eq:incompleteuonesample}) by taking $S$ as $S_{oa}=\{{\bm \eta}^1,\ldots,{\bm \eta}^m\}$. We now give a toy example of choosing $\eta^i$, for $i=1,\ldots,m$. Suppose $d=4$, $n=9$, and $$X_6\leq X_8\leq X_2\leq X_4\leq X_7\leq X_5\leq X_3\leq X_9\leq X_1.$$
Then, we have $L=3$ groups listed as
$G_1=\{6,8,2\},G_2=\{4,7,5\}$, and $G_3=\{3,9,1\}.$ An example of $OA(m=9,d=4,L=3,t=2)$ in step 1 is given as follows in transpose:
\begin{equation}\label{eqn:210}
A^T=
\left(
\begin{array}{ccccccccc}
1&1&1&{\bf 2}&2&2&3&3&3\\
1&2&3&{\bf 1}&2&3&1&2&3\\
1&2&3&{\bf 2}&3&1&3&1&2\\
1&2&3&{\bf 3}&1&2&2&3&1
\end{array}
\right).
\end{equation}
The fourth row of $A$, namely $(2,1,2,3)$, means we are sampling ${\bm \eta}^4$ from $G_2\times G_1\times G_2\times G_3$. One possible outcome for ${\bm \eta}^4$ could be $(4,8,7,3)$. Repeating this for each row of $A$, we could possibly have the $\mathcal{X}_{{\bm \eta}^i}$, for $i=1,\ldots,9$, used in the construction as follows:
\begin{equation}\label{eq:example}
\{\mathcal{X}_{{\bm \eta}^1},\ldots,\mathcal{X}_{{\bm \eta}^9}\}=\left\{
\begin{array}{ccccccccc}
X_6&X_8&X_2&X_{\bf 4}&X_4&X_5&X_9&X_3&X_1\\
X_2&X_4&X_3&X_{\bf 8}&X_7&X_9&X_8&X_5&X_9\\
X_8&X_6&X_9&X_{\bf 7}&X_1&X_2&X_1&X_2&X_4\\
X_6&X_5&X_1&X_{\bf 3}&X_6&X_4&X_7&X_9&X_6
\end{array}
\right\}.
\end{equation}

To proceed with the asymptotic properties of $U_{oa}$, we define
\begin{eqnarray}\label{eq:rt}
R(t)=\sum_{j>t}\binom{d}{j}\delta_{j}^2.
\end{eqnarray}
\begin{theorem}\label{thm:noassumptiononesample}
For any $(g,F)$, using $OA(m,d,L,t)$ in step 1 of the ICUDO algorithm, we have
\begin{eqnarray}\label{eq:noassumptiononesample}
{\rm MSE}(U_{oa})={\rm MSE}(U_0)+\frac{R(t)}{m}+o\left(\frac{1}{m}\right)+O\left(\frac{1}{n^2}\right).
\end{eqnarray}
\end{theorem}
We now explain the meanings of the three terms in (\ref{eq:noassumptiononesample}) generated in the process of approximating the complete U-statistic $U_0$ using $U_{oa}$. The term $O(n^{-2})$ is the bias square of $U_{oa}$ due to the inclusion of combinations with duplicate units, such as the first column of (\ref{eq:example}). Essentially, $U_{oa}$ is approximating the V-statistic, which is biased for $\Theta$ itself. The term $o(m^{-1})$ is due to the sampling variability when we draw one point from each selected group, that is, step 3 of the algorithm. The term $R(t)/m$ is due to the usage of the OA structure in place of a complete enumeration of all group combinations. Compared with the second term in (\ref{eq:generalchangeonesample}) for the ICUR, $R(0)/m$, we are able to eliminate all $\delta_j^2$ with $j\leq t$ owing to the projective uniformity of the OA in all $t$-dimensional projected spaces. If $\delta_j^2=0$ for $d'\leq j\leq d$, an OA with strength $t\geq d'$ yields $R(t)=0$. We discuss the hidden benefit of using a lower strength OA in Example \ref{example:tradeoff}. 

In the non-degenerate case, recall the MSE$(U_0)\asymp n^{-1}$ and $\lim_{n\rightarrow\infty}{\rm Eff}(U_{\rm RND})\leq d/(1+d)$ for the ICUR when $m\asymp n$. Under the same situation, Theorem \ref{thm:noassumptiononesample} implies that $U_{oa}$ is asymptotically efficient by simply taking $t=d$. In fact, stronger results can be derived for the ICUDO so that $m$ is allowed to grow more slowly than $n$ under various conditions. We give Theorem \ref{thm:lipschitzonesample} here as one example; additional results can be found in the Supplementary Material.
\begin{theorem}\label{thm:lipschitzonesample}
Suppose $(i)$ the kernel function $g$ is Lipschitz continuous, and $(ii)$ $F$ has density function $f(x)>c$ for some fixed $c>0$ and $x\in[a,b]$, and $f(x)=0$ otherwise. For $U_{oa}$ based on $OA(m,d,L,t)$ with $L^2\leq n(\log n)^{-1}$, we have
\begin{eqnarray}\label{eqn:1202}
{\rm MSE}(U_{oa})={\rm MSE}(U_0)+\frac{R(t)}{m}+O\left(\frac{1}{mL^2}\right)+O\left(\frac{1}{n^2}\right).
\end{eqnarray}
\end{theorem}

For $t=d=2$, we automatically have $R(t)=0$. If the conditions in Theorem \ref{thm:lipschitzonesample} hold, we only need $m\succ \sqrt{n}$ to achieve ${\rm Eff}(U_{oa})\rightarrow 1$, while the ICUR requires $m\succ n$. In general, $R(t)$ decreases in $t$ and could vanish if we take $t$ large enough so that $\delta_j^2=0$, for all $j>t$. Without knowledge of $\delta_j^2$, simply taking $t=d$ will eliminate $R(t)$ too. On the other hand, the term $O\left(\frac{1}{mL^2}\right)$ in (\ref{eqn:1202}) is decreasing in $L$, meaning the more groups we use to divide the data, the more homogeneous the units we could have in each group. However, $L$ and $t$ are subject to the constraint $m=\lambda L^t$, where $\lambda$ is the number of replicates of each $t$-tuple in OA and is equal to one in all examples presented here. As a result, $L$ and $t$ cannot be increased simultaneously. To gain insight to the trade-off between $L$ and $t$, we need to determine the constant term for $O\left(\frac{1}{mL^2}\right)$. For this, we derive the following theorem. A more detailed discussion on how to choose $L$ and $t$, given $m$, is provided in the Supplementary Material. Denote by $U(0,1)$ the uniform distribution on $[0,1]$.


\begin{theorem}\label{thm:tradeoff} 
Suppose $g$ has a continuous first-order derivative on $[0,1]^d$, $X\sim U(0,1)$, and there exists some $c\in (0,\frac{1}{2})$, such that $L\preceq n^{c}$. For $U_{oa}$ based on $OA(m,d,L,t)$,
\begin{eqnarray}\label{eq:tradeoff}
{\rm MSE}(U_{oa})={\rm MSE}(U_0)+\frac{R(t)}{m}+\frac{d}{12mL^2}E\gamma^2(X_1,\ldots,X_d)+o\left(\frac{1}{mL^2}\right),
\end{eqnarray}
where $\gamma(x_1,\ldots,x_d)=\frac{\partial g}{\partial x_1}(x_1,\ldots,x_d)$.
\end{theorem}

The assumption of a uniform distribution for $X$ is not as strict as it seems. To see this, for $X\sim F$, let $Z=F(X)\sim U(0,1)$. Applying Theorem \ref{thm:tradeoff} to $g_{F}(Z_1,\ldots,Z_d):=g(F^{-1}(Z_1),\ldots,F^{-1}(Z_d))=g(X_1,\ldots,X_d)$, we have the following corollary.
\begin{corollary}
Suppose $g_F$ has a continuous first-order derivative on $[0,1]^d$, and there exists some $c\in (0,\frac{1}{2})$, such that $L\preceq n^{c}$. Then, {\rm(\ref{eq:tradeoff})} still holds.
\end{corollary}

The term $E\gamma^2$ in (\ref{eq:tradeoff}) provides a nice interpretation of the trade-off between $t$ and $L$. When the kernel function $g$ has a large variability (large $E\gamma^2$), it is more challenging to make each group as homogeneous as possible, which enforces larger values of $L$. On the other hand, if $g$ is quite flat on the domain (small $E\gamma^2$), we prefer fewer groups to improve the strength of the OA.  
\begin{example}\label{example:tradeoff}
The kernel function $g(x_1,x_2,x_3)=x_1x_2x_3$ estimates $\mu^3$, where $\mu=E(X)$. We compare the performance of three methods: $U_{\rm RND}$; $U_{oa_2}$ based on $OA(m,3,\sqrt{m},2)$, with strength $t=2$; and $U_{oa_3}$ based on $OA(m,3,m^{1/3},3)$, with strength $t=3$. The data consist of $n=10^4$ i.i.d. observations simulated from 
$N(\mu,1)$, where $\mu$ takes the values of $0.5$ and $2$; see Table \ref{tb:t.normal} for the simulation results.
\end{example}

\begin{table}[h]
  \begin{center}
    \caption{Result of Example \ref{example:tradeoff}.}\label{tb:t.normal}
    \begin{tabular}{|c|ccc|ccc|}\hline
    
    \multirow{2}{*}{$m/n$}{}&
    \multicolumn{3}{c|}{$\mu=0.5$}&\multicolumn{3}{c|}{$\mu=2$}\\
    &Eff$(U_{\rm RND})$&${\rm Eff}(U_{oa_2})$&${\rm Eff}(U_{oa_3})$&Eff$(U_{\rm RND})$&${\rm Eff}(U_{oa_2})$&${\rm Eff}(U_{oa_3})$\\\hline
 0.005 &0.133\%&0.171\% &{\bf0.218}\%& 1.110\% &{\bf9.908}\%&2.323\%  \\
  0.01 &0.290\%&{0.464}\%&{\bf0.579}\%& 2.485\% &{\bf26.84}\%&8.455\%  \\
  0.05 &1.291\%&{2.448}\%&{\bf6.096}\%&10.31\%&{\bf75.12}\% &51.71\%  \\
  0.1 &2.936\%&{4.527}\%&{\bf16.62}\%&20.13\% &{\bf91.87}\%&76.80\%  \\
  0.5 &12.58\%&{21.89}\%&{\bf71.78}\%&50.78\%&{\bf100.0}\% &98.53\% \\
  1.0 &21.05\%&{33.26}\%&{\bf99.94}\%&67.51\% &{\bf100.0}\%&99.64\%\\\hline

    \end{tabular}
  \end{center}
\end{table}

In Table \ref{tb:t.normal}, both $U_{oa_2}$ and $U_{oa_3}$ outperform $U_{\rm RND}$ significantly. The advantage of the ICUDO over the ICUR is discussed below in additional examples. Furthermore, we find that the winning strategy changes from $U_{oa_3}$ to $U_{oa_2}$ as we increase the mean $\mu$ of the distribution. This observation well illustrates the comments after Theorem \ref{thm:tradeoff} on the relevance of $E\gamma^2$ in determining the optimal value of the strength $t$. That is, for larger $E\gamma^2$, we are more inclined to choose a smaller strength. This is validated by our second observation together with $E\gamma^2=(\mu^2+1)^2$, which increases in $\mu(>0)$. 

Note that the applicability of Theorem \ref{thm:lipschitzonesample} and its variants, Theorems \ref{thm:clipschitzonesample}--\ref{thm:linearcombination} in the Supplementary Material is broader than it appears. To see this, let $\phi:R\rightarrow R$ be a one-to-one mapping. Denote by $F_{\phi}$ the distribution of the transformed random variable $Z=\phi(X)$, which leads to the following representation: 
$$g_{\phi}(z_1,\ldots,z_d):=g(\phi^{-1}(z_1),\ldots,\phi^{-1}(z_d))=g(x_1,\ldots,x_d).$$
If ($g_{\phi},F_{\phi}$) satisfies the conditions in these theorems, corresponding results also hold for the pair ($g,F$). For example, suppose $g(x_1,x_2)=x_1^{-a}x_2^{-a}$ and $F$ is a Pareto distribution with shape and scale parameters $a$ and $b$, respectively. The Pareto distribution is neither light-tailed nor bounded, and hence violates the conditions in Theorem \ref{thm:lipschitzonesample}. By taking $\phi(x)=1-(b/x)^a$, we have $\phi(X)\sim U(0,1)$. It can be verified that the conditions in Theorem \ref{thm:lipschitzonesample} are satisfied by ($g_{\phi},F_{\phi}$).

\noindent{\bf 2.2. Multi-sample U-statistics}

For $k=1,\ldots,K$, let
$X_{1}^{(k)},\ldots,X_{n_k}^{(k)}$ be a random sample of size $n_k$ from the distribution $F_k$. The UMVUE of 
$$\Theta=\int g(x^{(1)}_1,\ldots,x^{(1)}_{d_1},\cdots,x^{(K)}_1,\ldots,x^{(K)}_{d_K})dF_1(x^{(1)}_1)\ldots dF_K(x^{(K)}_{d_K})$$
is given by the generalized U-statistic
$$U_0=\prod_{k=1}^K\binom{n_k}{d_k}^{-1}\sum_{\bm\eta\in\prod_{k=1}^KS_{n_k,d_k}}g(\mathcal{X}_{\bm\eta}),$$
$$S_{n_k,d_k}=\{{\bm \eta}_k=(\eta_{k,1},\ldots,\eta_{k,d_k}): 1\leq \eta_{k,1}<\eta_{k,2}<\ldots<\eta_{k,d_k}\leq n_k\},$$
$$\mathcal{X}_{{\bm\eta}}=(\mathcal{X}_{\bm\eta_1},\ldots,\mathcal{X}_{\bm\eta_K})=(X^{(1)}_{{ \eta}_{1,1}},\ldots,X^{(1)}_{{ \eta}_{1,d_1}},\cdots,X^{(K)}_{{ \eta}_{K,1}},\ldots,X^{(K)}_{{ \eta}_{K,d_K}}).$$
The $d(=\sum_{k=1}^Kd_k)$-dimensional kernel function $g$ is symmetric about any $d_k$-dimensional sub-input $\{x^{(k)}_1,\ldots,x^{(k)}_{d_k}\}$. The generalized U-statistic reduces to the traditional U-statistic when $K=1$.
An incomplete generalized U-statistic is given by
\begin{equation}\label{eq:incompleteu}
U=\frac{1}{m}\sum_{{\bm\eta}\in S}g(\mathcal{X}_{\bm\eta}),
\end{equation}
where $S\subset\prod_{k=1}^KS_{n_k,d_k}$ and $m=|S|$. We construct the multi-sample ICUDO as follows. For ease of illustration, we assume $n_k$' is a multiple of $L$.

\begin{itemize}
\item[Step 1.] Let $A_0$ be an $OA(m,d,L,t)$. Adopt random level permutations $\{\pi_1,\ldots,\pi_d\}$ of columns of $A_0$ independently. Specifically, for each $l\in\mathcal{Z}_L$, change all elements $l$ in the $j$th column of $A_0$ to $\pi_j(l)$. The $m$ rows of the resulting array $A$ are denoted by $\{\bm a^i=(\bm a^{i}_1,\ldots,\bm a^{i}_K):i=1,\ldots,m;\bm a^{i}_k\in\mathcal{Z}_L^{d_k},k=1,\ldots,K\}$.
\item[Step 2.] For each $k=1,\ldots,K$, create the partition $\mathcal{Z}_{n_k}=\bigcup_{l=1}^LG^{(k)}_l$, such that $|G^{(k)}_l|=n_kL^{-1}$ for $l\in\mathcal{Z}_L$, and $X_{i_1}^{(k)}\leq X_{i_2}^{(k)}$ for any $i_1\in G_{l_1}^{(k)}$, $i_2\in G_{l_2}^{(k)}$, with $l_1<l_2$. For any ${\bm a}=({\bm a}_1,\ldots,{\bm a}_K)$ with ${\bm a}_k=(a_{k,1},\ldots,a_{k,d_k})\in\mathcal{Z}_L^{d_k}$, define
\begin{equation}\label{eq:Gv}
\mathcal{G}_{{\bm a}}=\prod_{k=1}^K\prod_{j=1}^{d_k}G^{(k)}_{a_{k,j}}.
\end{equation}
\item[Step 3.] For $i=1,\ldots,m$, independently draw an element ${\bm \eta}^i$ uniformly from $\mathcal{G}_{{\bm a}^i}$, where $\bm a^i$ is the $i$th row of $A$:
\begin{equation}\label{uoamulti}
U_{oa}=\frac{1}{m}\sum_{i=1}^mg(\mathcal{X}_{\bm\eta^i}).
\end{equation}
\end{itemize}
An example is given in the Supplementary Material.
For any $j_{k,1},\ldots,j_{k,d_k}\in \mathcal{Z}_{d_k}$ and $k\in\mathcal{Z}_K$, assume
$$Eg^2\left(X^{(1)}_{j_{1,1}},\ldots,X^{(1)}_{j_{1,d_1}},\cdots,X^{(K)}_{j_{K,1}},\ldots,X^{(K)}_{j_{K,d_K}}\right)<\infty.$$
Let $n_{\min}=\min\{n_1,\ldots,n_K\}$ and $n_{\max}=\max\{n_1,\ldots,n_K\}$.
Here, we assume $n_{\min}\asymp n_{\max}$ and $L\prec n_{\min}$.
Let ${\bm u}=({\bm u}_1,\ldots,{\bm u}_K)$, where ${\bm u}_k\subseteq \mathcal{Z}_{d_k}$. Define $dF_{{\bm u}}=\prod_{k=1}^K\prod_{j\in{\bm u}_k}dF_k(x_{j}^{(k)})$. For any ${\bm u}$ and ${\bm x}=(x^{(1)}_1,\ldots,x^{(1)}_{d_1},\cdots,x^{(K)}_1,\ldots,$ $x^{(K)}_{d_K})$, we recursively define
\begin{equation}\label{eq:gufunction}
g_{\bm u}({\bm x})=\int g({\bm x})dF_{{\bm u}^c}~~~~h_{\bm u}({\bm x})=g({\bm x})-\sum_{{\bm v}\subset{\bm u}}h_{\bm v}(x),\notag
\end{equation}
where $\bm u^c=(\bm u_1^c,\ldots,\bm u_K^c)=(\mathcal{Z}_{d_1}\setminus\bm u_1,\ldots,\mathcal{Z}_{d_K}\setminus\bm u_K)$, $g_{\emptyset}({\bm x})=\Theta$ and $h_{\emptyset}({\bm x})=0$, ${\bm v}=({\bm v}_1,\ldots,{\bm v}_K)$, and ${\bm v}\subset{\bm u}$ means ${\bm v}_k\subseteq {\bm u}_k$ (${\bm v}\neq{\bm u}$).

For $\bm u$, we can define
$\sigma_{{\bm u}}^2={\rm Var}(g_{\bm u})~{\rm and}~\delta_{{\bm u}}^2={\rm Var}(h_{\bm u})$.
The MSE of the complete generalized U-statistic is given by \cite{sen:1974} as
\begin{equation}\label{eq:generalvar}
{\rm MSE}(U_0)=\prod_{k=1}^K\binom{n_k}{d_k}^{-1}\sum_{\bm u=(\bm u_1,\ldots,\bm u_K)}\left\{\prod_{k=1}^K\binom{d_k}{|\bm u_k|}\binom{n_k-d_k}{d_k-|\bm u_k|}\right\}\sigma_{\bm u}^2.\notag
\end{equation}
Let $|\bm u|=\sum_{k=1}^K|\bm u_k|$. The generalized U-statistic and the kernel function are called {\it order-q degenerate} if $\sigma_{\bm u}^2=\sum_{\bm v\in\bm u}\delta_{\bm v}^2=0$, for all $|\bm u|\leq q$, and there exists $\bm u'$ such that $\sigma_{\bm u'}^2>0$ and $|\bm u'|=q+1$. We have ${\rm MSE}(U_0)=O(n^{-(q+1)})$ in this case. For the non-degenerate case $q=0$, we have ${\rm MSE}(U_0)\asymp n^{-1}$ . With a slight abuse of notation, let $\sigma_{(j_1,\ldots,j_K)}=\sigma_{\bm u}$ and $\delta_{(j_1,\ldots,j_K)}=\delta_{\bm u}$, for $\bm u=(\bm u_1,\ldots,\bm u_K)$, with $|\bm u_k|=j_k$, $k=1,\ldots,K$. For the ICUR, we have
\begin{equation}\label{eq:generalchange}
{\rm MSE}(U_{\rm RND})={\rm MSE}(U_0)+\frac{R(0)}{m}+O\left(\frac{1}{mn_{\min}}\right),\notag
\end{equation}
$$R(t)=\sum_{\bm u:|\bm u|>t}\delta_{\bm u}^2=\sum_{j_1=0}^{d_1}\cdots\sum_{j_K=0}^{d_K}I(j_1+\cdots+j_K>t)\prod_{k=1}^K\binom{d_k}{j_k}\delta_{(j_1,\ldots,j_K)}^2.$$ 
The last term above reduces to the form of $R(t)$ for the one-sample case, but the second term yields a parsimonious presentation for the multi-sample case. The corresponding properties of $U_{oa}$ are given as follows.
\begin{theorem}\label{thm:noassumption}
For $U_{oa}$ based on $OA(m,d,L,t)$, for any pair of $(g,F)$, we have
\begin{eqnarray}\label{eq:noassumptionmulti}
{\rm MSE}(U_{oa})={\rm MSE}(U_0)+\frac{R(t)}{m}+o\left(\frac{1}{m}\right)+O\left(\frac{1}{n_{\min}^2}\right).
\end{eqnarray}
\end{theorem}

Theorem \ref{thm:noassumption} is basically a multi-sample version of Theorem \ref{thm:noassumptiononesample}, and its result can be strengthened in the same way. The details are omitted here to conserve space. We conclude this section with a machine learning example.

\begin{figure}
\centering
\begin{minipage}[t]{0.49\textwidth}
\centering
\includegraphics[width=6.7cm]{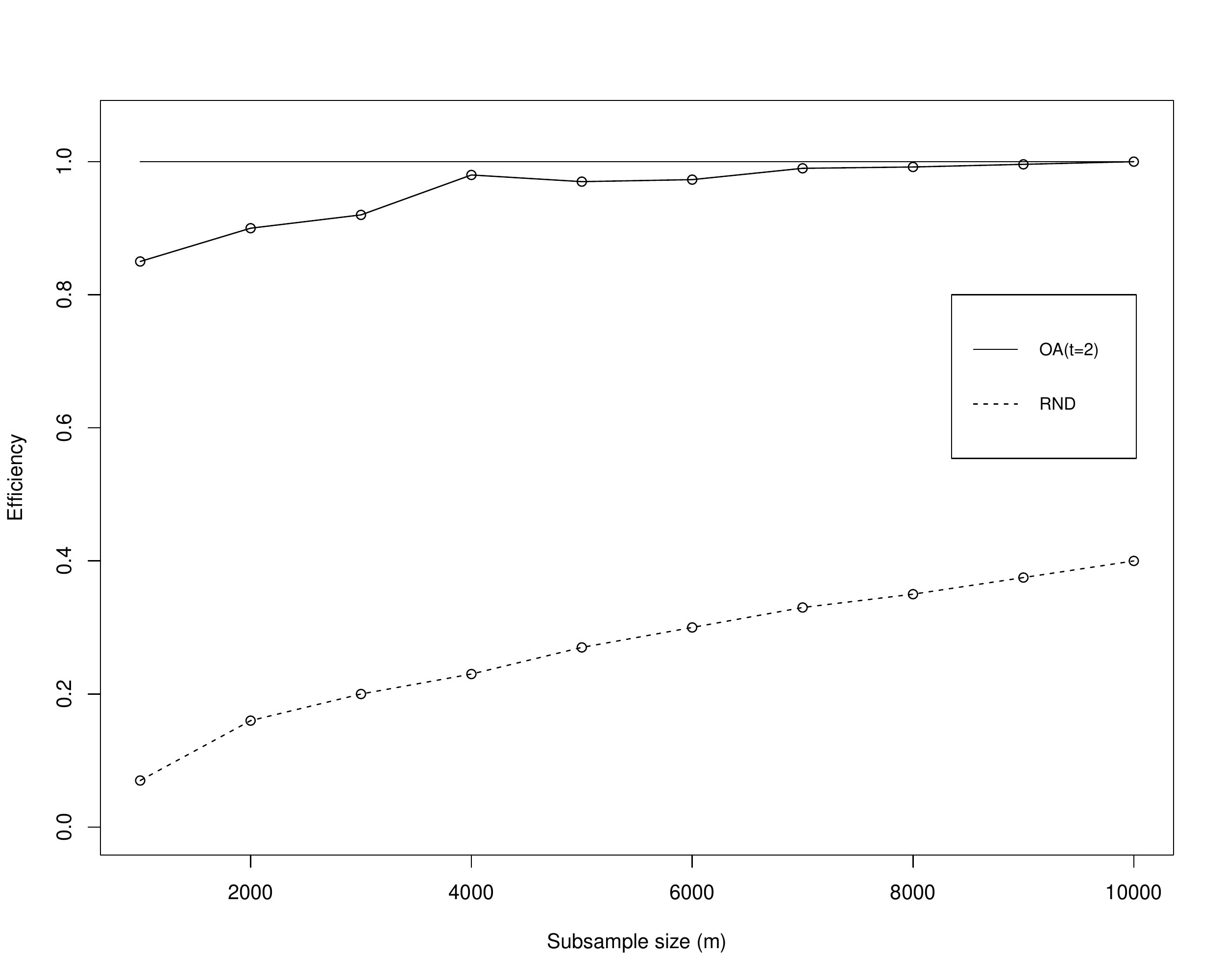}
\end{minipage}
\begin{minipage}[t]{0.49\textwidth}
\centering
\includegraphics[width=6.7cm]{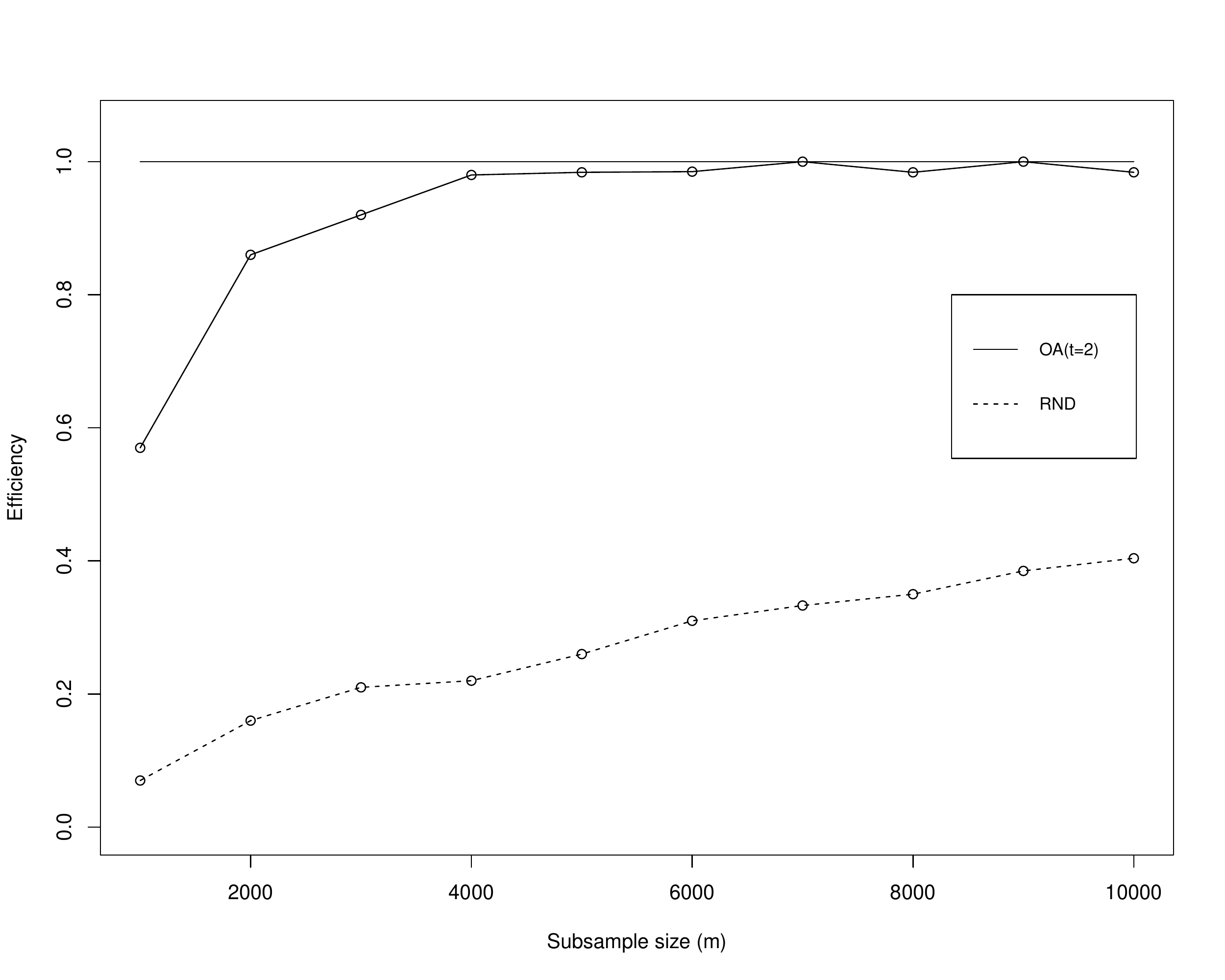}
\end{minipage}
\caption{Comparison of efficiencies of $\tilde{U}_{oa}$ and $U_{\rm RND}$ with respect to subsample size $m$ for loss function ($i$) (left) and ($ii$) (right).}
\label{fig:multid}
\end{figure}

\begin{example} \label{example:rankingloss}(Ranking measure, \cite{chen:2009}). The ranking measure is an important topic in machine learning research. In the commonly used pairwise approach, the loss for a given classifier score function $f$ is given by
$$L(f)=\sum_{1\leq i<j\leq K}\sum_{x\in G_i,y\in G_j}\psi(f(y)-f(x)),$$
where $G_1,\ldots,G_K$ are $K$ groups ranked in ascending order. Here, $\psi$ could that the form of
\begin{itemize} 
\item[$(i)$] hinge function: {\rm$\psi(z)=(1-z)_{+}$}, or a
\item[$(ii)$] logistic function: {\rm$\psi(z)=\log(1+\exp(-z))$} 
\end{itemize}
for the Ranking SVM and RankNet methods, respectively. In the simulation, we set $K=2$, that is, the two-sample case, $|G_1|=|G_2|=10^4$, $f(G_1)\sim N(0,4)$, and $f(G_2)\sim N(5,4)$.
Figure \ref{fig:multid} reveals the high efficiency of $\tilde{U}_{oa}$ compared with that of $U_{\rm RND}$.
\end{example}

\vspace{0.5cm}
\setcounter{section}{3} 
\setcounter{equation}{2} 

\lhead[\footnotesize\thepage\fancyplain{}\leftmark]{}\rhead[]{\fancyplain{}\rightmark\footnotesize\thepage}
\noindent{\large\bf 3. Debiased ICUDO for degenerate cases}

Recall the ICUDO procedure is actually biased owing to the inclusion of combinations with duplicate units. The bias square is $O(n^{-2})$ for any pair $(g,F)$, which is negligible compared to ${\rm Var}(U_0)\asymp n^{-1}$ in the non-degenerate case. One can see that it is no longer negligible in the degenerate case. In this section, we propose a debiased version of the ICUDO. 

We provide details for the multi-sample cases, where the one-sample cases are achieved by taking $K=1$. To proceed, Let $S_0^*=\{({\bm \eta}_1,\ldots,{\bm \eta}_K):  {\bm \eta}_k=(\eta_{k,1},\ldots,\eta_{k,d_k})\in\mathcal{Z}_{n_k}^{d_k}, \eta_{k,j_1}\neq \eta_{k,j_2} {\rm~for~any~} j_1\neq j_2 \}$. The debiased ICUDO is constructed in the same way as the original, except that step 3 changes as follows:
\begin{itemize}
\item[Step 3$'$.] For $i=1,\ldots,m$, independently draw ${\bm \eta}^i$ from the uniform distribution on $\mathcal{G}_{{\bm a}^i}\cap S_0^*$.
Adopting (\ref{eq:incompleteu}) with $S_{oa}^*=\{{\bm \eta}^1,\ldots,{\bm \eta}^m\}$, we have the debiased ICUDO as
\begin{equation}\label{uoade}
\tilde{U}_{oa}=\frac{1}{m}\sum_{i=1}^m\omega_{\bm{\eta}^i}g(X_{\bm \eta^{i}}),
\end{equation}
where $\omega_{\bm \eta^i}=L^d|\mathcal{G}_{{\bm a}^i}\cap S_0^*|/|S_0^*|.$
\end{itemize}

\begin{theorem}\label{thm:debias}
$\tilde{U}_{oa}$ based on $OA(m,d,L,t)$ is an unbiased estimator, and
\begin{eqnarray}\label{eq:debiased}
{\rm MSE}(\tilde{U}_{oa})={\rm MSE}(U_0)+\frac{R(t)}{m}+o\left(\frac{1}{m}\right).
\end{eqnarray}
\end{theorem}
Theorem \ref{thm:debias} is analogous to Theorems \ref{thm:noassumptiononesample} and \ref{thm:noassumption} for the one-sample and multi-sample cases, respectively, except that the bias square term $O(n^{-2})$ and $O(n_{\min}^{-2})$ are eliminated. Now, for an order-$q$ degenerate U-statistic, the debiased ICUDO can be asymptotically efficient with $m\asymp n^{q+1}$, while the ICUR requires $m\succ n^{q+1}$. Moreover, we could allow $m$ to grow more slowly for the debiased ICUDO under some mild conditions on $(g,F)$. For example, when $d=2$, $q=1$, and the conditions of Theorem \ref{thm:lipschitzonesample} hold, the debiased ICUDO only needs $m\succ n$ to be asymptotically efficient, while the ICUR requires $m\succ n^2$. For the general order $q$ of degeneration, we have $m^*_{oa}=(m^*_{\rm RND})^{\frac{d}{d+1}}$, for all $d$, under the conditions in Theorem \ref{thm:lipschitzonesample}. Here, $m^*_{oa}$ and $m^*_{\rm RND}$ represent the minimum $m$ required for the ICUDO and ICUR, respectively, to be asymptotically efficient. 

We conclude this section with the following multi-sample example. The kernel function is degenerate, and hence favors a debiased ICUDO. However, the highest order $\delta^2$-value vanishes, which encourages a lower strength of OA. The comparison is made between the ICUR and different versions of the ICUDO.

\begin{example}\label{example:final}
Let $K=2$, $d_1=d_2=2$, $d=4$, and 
$$g(x^{(1)}_1,x^{(1)}_2,x^{(2)}_1,x^{(2)}_2)=I(x^{(1)}_1<x^{(2)}_1,x^{(1)}_2<x^{(2)}_1)+I(x^{(2)}_1<x^{(1)}_1,x^{(2)}_2<x^{(1)}_1).$$
The construction of $U_{oa}$ and the debiased $\tilde{U}_{oa}$ is based on $OA(m,4,m^{1/3},3)$ and $OA(m,4,$ $m^{1/4},4)$.
For continuous distributions $F_1$ and $F_2$, it can be verified that
$$E g(X^{(1)}_1,X^{(1)}_2,X^{(2)}_1,X^{(2)}_2)=\frac{2}{3}+\int(F_1(x)-F_2(x))^2d(F_1(x)+F_2(x))/2,$$
which indicates the similarity of $F_1$ and $F_2$. The null hypothesis of $F_1=F_2$ is rejected when the U-statistic is significantly larger than $2/3$. Note that the corresponding U-statistic is degenerate under the null hypothesis. See Table \ref{tab:final} for the simulation results when both samples are simulated from $N(0,1)$ with sample sizes $n_1=n_2=10^3$.
\end{example}

Note that in the $g$ function of Example \ref{example:final}, the two separate parts are all functions of three inputs. Thus, $R(4)=0$, and we can claim that $t=3$ works better than $t=4$, which is verified by the results in Table \ref{tab:final}.

\begin{table}[h]
    
  \begin{center}\tabcolsep 6.6pt
    \caption{Result of Example \ref{example:final}.}\label{tab:final}\tabcolsep 5pt
    \begin{tabular}{|c|cccccccc|}\hline
   $m/\binom{n}{2}$ &0.002   &0.01   &0.02  &0.04 &0.06 &0.1 &0.14 &0.2 \\\hline

   ${\rm Eff}(\tilde{U}_{oa_3})$ &0.836\% &{10.9\%} &{\bf 15.6\%} &{\bf 35.9\%} &{\bf 44.9\%} &{\bf 56.9\%} &{\bf 75.1\%} &{\bf 94.1\%} \\
    \hline
   ${\rm Eff}(U_{oa_3})$  &{0.861\%} &9.50\% &12.9\% &25.2\% &28.3\% &29.8\% &36.3\% &39.0\% \\
    \hline
   ${\rm Eff}(U_{oa_4})$  &{ 0.450\%} &4.96\% &6.78\% &10.6\% &10.7\% &11.9\% &14.5\% &15.6\% \\
    \hline
   ${\rm Eff}(U_{\rm RND})$  &0.179\% &0.701\% &1.50\% &2.93\% &4.19\% &7.84\% &10.9\% &13.1\% \\

   \hline

    \end{tabular}
  \end{center}

\end{table}

\setcounter{section}{4} 
\setcounter{equation}{3} 

\lhead[\footnotesize\thepage\fancyplain{}\leftmark]{}\rhead[]{\fancyplain{}\rightmark\footnotesize\thepage}
\noindent{\large\bf 4. ICUDO for multi-dimensional data}

Note that step 2 of the ICUDO algorithm in Section 2 does not apply to multi-dimensional data because it relies on ordering the univariate data. To remedy this, we adopt a clustering algorithm to divide the data into homogeneous groups. In this regard, the clustered group sizes may vary. This will necessitate a re-weighting procedure similar to the debiasing step in Section 3. To save space, we focus on the debiased ICUDO and adopt the notation of the multi-sample U-statistics in the study of multi-dimensional data. For $k=1,\ldots,K$, let $X_{1}^{(k)},\ldots,X_{n_k}^{(k)}$ be a random sample of size $n_k$ from the multi-dimensional distribution $F_k$. The algorithm is given as follows.

\begin{itemize}
\item[Step 1.] Let $A_0$ be an $OA(m,d,L,t)$. Adopt random level permutations $\{\pi_1,\ldots,$ $\pi_d\}$ of columns of $A_0$ independently. Specifically, for $l\in\mathcal{Z}_L$, change all elements $l$ in the $j$th column of $A_0$ to $\pi_j(l)$. The $m$ rows of the resulting array $A$ are denoted by $\{\bm a^i=(\bm a^{i}_1,\ldots,\bm a^{i}_K):i=1,\ldots,m;\bm a^{i}_k\in\mathcal{Z}_L^{d_k},k=1,\ldots,K\}$.
\item[Step 2.] Let $\mathcal{P}^{(k)}=\{G_1^{(k)},\ldots,G_L^{(k)}\}$ denote an $L$-group partition from the clustering of $\{X_{1}^{(k)},\ldots,X_{n_k}^{(k)}\}$. For any ${\bm a}=({\bm a}_1,\ldots,{\bm a}_K)$, with ${\bm a}_k=(a_{k,1},\ldots,a_{k,d_k})\in\mathcal{Z}_L^{d_k}$, define
\begin{equation}\label{eq:Gv}
\mathcal{G}_{{\bm a}}=\prod_{k=1}^K\prod_{j=1}^{d_k}G^{(k)}_{a_{k,j}}.
\end{equation}
\item[Step 3.] For $i=1,\ldots,m$, independently draw an element ${\bm \eta}^i$ uniformly from $\mathcal{G}_{{\bm a}^i}$, where $\bm a^i$ is the $i$th row of $A$. Let $\omega_{\bm \eta^i}=L^d|\mathcal{G}_{{\bm a}^i}\cap S_0^*|/|S_0^*|$.
\begin{equation}\label{uoamulti}
\tilde{U}_{oa}=\frac{1}{m}\sum_{i=1}^m\omega_{\bm \eta^i}g(X_{\bm\eta^i}).
\end{equation}
\end{itemize}
An example of the construction is given in the Supplementary Material.
\begin{theorem}\label{thm:multid1} Suppose $\omega_{\bm \eta^i}\rightarrow 1$ uniformly as $n,L\rightarrow\infty$.
For $\tilde{U}_{oa}$ based on $OA(m,d,L,t)$, we have
\begin{eqnarray}\label{eq:multid}
{\rm MSE}(\tilde{U}_{oa})={\rm MSE}(U_0)+\frac{R(t)}{m}+o\left(\frac{1}{m}\right).
\end{eqnarray}
\end{theorem}
The $R(t)$ in (\ref{eq:multid}) is given by (\ref{eq:rt}), except that the univariate distribution $F$ is changed to a multi-dimensional distribution. The assumption in Theorem \ref{thm:multid1} naturally holds if we force balance the group size in the clustering process. By applying the full strength $t=d$ OA to Theorem \ref{thm:multid1}, we have the following corollary.

\begin{corollary}\label{thm:multid2} 
For $\tilde{U}_{oa}$ based on $OA(m,d,L,d)$, for any pair of $(g, F)$, we have
\begin{eqnarray}\label{eq:noassumptionmultid}
{\rm MSE}(\tilde{U}_{oa})={\rm MSE}(U_0)+o(m^{-1}).
\end{eqnarray}
\end{corollary}


The choice of $t$ has been discussed and is illustrated in Examples \ref{example:tradeoff} and \ref{example:final}. We do not compare different $t$ in the following examples because $d=2$ always holds, and so $t\leq 2$. We always take $t=2$, $L=10,20,\ldots,100$, and $m=L^t$.

\begin{example} \label{example:kendall}(Kendall's tau, \cite{chen:2018}). The Kernel function $h((x_1, y_1),(x_2, $ $y_2)) = 2 I(x_1 < x_2, y_1 < y_2) + 2 I(x_2 <
x_1, y_2 < y_1)-1$. For simplicity, we assume that $(X,Y)$ follows a normal distribution, with $\mu=(0,0)$ and $\Sigma={\rm diag}(3,1)$. Set $n=10^4$. The MSE when estimating the Kendall correlation using $U_{\rm RND}$ and $\tilde{U}_{oa}$ is shown in Table \ref{tab:kendall}. As a reference, we have ${\rm MSE}(U_0)=8.97\times 10^{-5}$.
\end{example}

\begin{table}[h]
    
  \begin{center}\tabcolsep 4.5pt
    \caption{Result of Example \ref{example:kendall}.}\label{tab:kendall}\tabcolsep 3.5pt
    \begin{tabular}{|c|cccccccccc|}\hline
   $m$ &100   &400   &900  &1600 &2500 &3600 &4900 &6400 &8100 &10000 \\\hline

   ${\rm MSE}({U}_{\rm RND})$ &.765 &.191 &.0903 &.0515 &.0260 &.0195 &.0167 &.0137 &.0098 &.0089 \\
    \hline
   ${\rm MSE}(\tilde{U}_{oa})$  &.075 &.0096 &.0032 &.0015 &.00063 &.00035 &.00023 &.00014 &.00011 &.00009 \\

   \hline

    \end{tabular}
  \end{center}

\end{table}

\begin{example} \label{example:monotonicity}(Testing stochastic monotonicity, \cite{lee:2009}). Let $(X,Y)$ be a real-valued random vector, and denote by $F_{Y|X}(y|x)$ the conditional distribution function of $Y$, given $X$. Consider the problem of testing the stochastic monotonicity hypothesis
$$H_0: F_{Y|X}(y|x)\leq F_{Y|X}(y|x'), \forall y\in R ~{\rm and ~whenever}~ x\geq x'.$$
This essentially tests where an increase in $X$ would induce an increase in $Y$ (e.g., income vs. expenditure in a household). \cite{lee:2009} proposed the following testing statistic:
\begin{eqnarray}\label{eqn:309}
U_n(x,x')=\frac{1}{n(n-1)}\sum_{1\leq i\neq j\leq n}(I\{Y_i\leq x'\}-I\{Y_j\leq x'\}){\rm sign}(X_i-X_j)K(x-X_i)K(x-X_j),~
\end{eqnarray}
where $K(x)=0.75(1-x^2)$. We simulate $(X,Y)$ from a normal distribution with $\mu=(0,0)$ and $\Sigma={\rm diag}(3,1)$, and calculate (\ref{eqn:309}) at $(x,x')=(0,0)$. For $n=10^4$, the comparison between $\tilde{U}_{oa}$ and $U_{\rm RND}$ is given in Table \ref{tab:mono}. As a reference, we have ${\rm MSE}(U_0)=2.572$.
\end{example}

\begin{table}[h]
    
  \begin{center}\tabcolsep 5.8pt
    \caption{Result of Example \ref{example:monotonicity}.}\label{tab:mono}\tabcolsep 4.7pt
    \begin{tabular}{|c|cccccccccc|}\hline
   $m$ &100   &400   &900  &1600 &2500 &3600 &4900 &6400 &8100 &10000 \\\hline

   ${\rm MSE}({U}_{\rm RND})$ &302.7 &69.01 &38.01 &17.45 &12.86 &8.613 &7.438 &6.273 &4.886 &4.327 \\
    \hline
   ${\rm MSE}(\tilde{U}_{oa})$  &33.18 &15.73 &8.848 &4.252 &3.524 &3.168 &2.732 &2.662 &2.630 &2.602 \\

   \hline

    \end{tabular}
  \end{center}

\end{table}

{
\begin{example} \label{example:clusterperformance}(Clustering performance evaluation, \cite{papa:2015}). For a given distance $D:\mathcal{X}\times\mathcal{X}\rightarrow R$ defined on $\mathcal{X}$, the performance of a partition ${P}$ can be evaluated from the data $X_1,\ldots,X_n\in\mathcal{X}$ using 
\begin{eqnarray}\label{eqn:311}
W({P})=\sum_{1\leq i<j\leq n}D(X_i,X_j)\cdot\sum_{\mathcal{C}\in{P}}I\{(X_i,X_j)\in\mathcal{C}^2\}.
\end{eqnarray}
Our purpose is to compare the different incomplete U-statistics of (\ref{eqn:311}); here, we focus on the k-means method for the comparison. The data are generated from a normal distribution with $\mu=(0,0)$ and $\Sigma={\rm diag}(1,2)$, and we divide the data into two groups. The MSE of $U_{\rm RND}$ and $\tilde{U}_{oa}$ when estimating $W(P)$ for different $m$ is shown in Table \ref{tab:clus}. As a reference, we have ${\rm MSE}(U_0)=1.043\times 10^{-4}$.
\end{example}

\begin{table}[h]
    
  \begin{center}\tabcolsep 4.5pt
    \caption{Result of Example \ref{example:clusterperformance}.}\label{tab:clus}\tabcolsep 3.7pt
    \begin{tabular}{|c|cccccccccc|}\hline
   $m$ &100   &400   &900  &1600 &2500 &3600 &4900 &6400 &8100 &10000 \\\hline

   ${\rm MSE}({U}_{\rm RND})$ &.216 &.0625 &.0346 &.0171 &.0064 &.0047 &.0038 &.0021 &.0017 &.0010 \\
    \hline
   ${\rm MSE}(\tilde{U}_{oa})$  &.011 &.0064 &.0038 &.0019 &.00056 &.00051 &.00038 &.00027 &.00013 &.00012 \\

   \hline

    \end{tabular}
  \end{center}

\end{table}}

\setcounter{section}{5} 
\setcounter{equation}{4} 

\lhead[\footnotesize\thepage\fancyplain{}\leftmark]{}\rhead[]{\fancyplain{}\rightmark\footnotesize\thepage}

\vspace{1cm}
\noindent{\large\bf 5. Conclusion}
\label{sec:discuz}

To tackle the computational issue of U-statistics, we have introduced a new type of incomplete U-statistic called the ICUDO, which has much higher efficiency than existing methods. The required computational burden, as indexed by the number of combinations $m$ for the ICUDO to be statistically equivalent to the complete U-statistic, is of smaller magnitude than existing methods. This was validated theoretically and empirically for degenerate and non-degenerate one- and multi-sample U-statistics on univariate and multi-dimensional data. In fact, $m$ is allowed to grow more slowly than the data size $n$ in the non-degenerate case.  

The OA plays a critical role in the construction of the ICUDO, in light of its projective uniformity. Other space-filling design schemes exist with similar properties, such as the OA-based Latin hypercube by \cite{tang:1993}, and the strong orthogonal array by \cite{he:2012}, which is used frequently in the design of computer experiments. By exhaustive simulations, we find the improvement of the efficiency by these design schemes over that of the ICUDO to be within 1\%. However, this improvement is not sufficient to advocate using these structures, owing to the extra complexity of the computation. Other improvements over the OA are based on optimal criteria, such as the generalized minimum aberration OA. However, no theoretical results are available for these fixed structures.

Lastly, the following offer potential future research directions. ($i$) For high-dimensional data, dimension-reduction techniques need to be integrated into our current algorithm. ($ii$) For multi-sample cases, we may divide different samples into different numbers of groups in some optimal way. This will induce more complicated OA structures. ($iii$) For the purpose of statistical inference, it would be of interest to study the asymptotic distributions of the ICUDO under different conditions. ($iv$) The dimension of the kernel functions is fixed at $d$ as $n$ increases, and all data are generated independently. In one important type of U-statistic based on stochastic processes, $d$ increases with $n$ and the data can be dependent. These topics will involve quite different methodologies, and hence are left to future work.

\setcounter{section}{6} 
\setcounter{equation}{4} 

\lhead[\footnotesize\thepage\fancyplain{}\leftmark]{}\rhead[]{\fancyplain{}\rightmark\footnotesize\thepage}

\vspace{1cm}
\noindent{\large\bf Supplementary Material}
\label{sec:supplement}

The online Supplementary Material generalizes the result of Theorem \ref{thm:lipschitzonesample} under additional conditions. It also provides details on how to choose the combination of $L$ and $t$ and illustrates the generation of the ICUDO for multi-sample and multi-dimensional cases.

\par
\vskip 14pt
\noindent {\large\bf Acknowledgments}

Dr. Kong's research was partially supported by NSFC grant 11801033 and the Beijing Institute of Technology Research Fund Program for Young Scholars.
Dr. Zheng's research was partially supported by the National Science Foundation, DMS-1830864.
\par
\vspace{.5cm}
\setcounter{section}{6} 
\setcounter{equation}{5} 

\lhead[\footnotesize\thepage\fancyplain{}\leftmark]{}\rhead[]{\fancyplain{}\rightmark\footnotesize\thepage}
\noindent{\large\bf Appendix. Proof of Theorems}
\label{sec:discuz}

Lemmas \ref{lem:barfunctionproperty}--\ref{lem:lusin} contribute to the proof of Theorem \ref{thm:noassumptiononesample}. Theorems \ref{thm:noassumption} and \ref{thm:multid1} can be proved similarly as Theorem \ref{thm:noassumptiononesample}, but only with more tedious analysis, and hence they are omitted due to the limit of space. For any $\bm a\in \mathcal{Z}_L^d$, we call the set $\mathcal{G}_{\bm a}=\prod_{j=1}^{d}G_{a_{j}}$ a {\it grid}. Let $F_n$ be the empirical distribution of $\{X_1,\ldots,X_n\}$ and define
$V=\int g(x_1,\ldots,x_d)dF_n(x_1)\ldots $ $dF_n(x_d).$
For given $F_n$ and ${\bm \eta}\in \mathcal{G}_{\bm a}$, define 
$$\bar{g}(\mathcal{X}_{\bm \eta})=|\mathcal{G}_{\bm a}|^{-1}\sum_{\bm{\eta}'\in \mathcal{G}_{\bm a}}g(\mathcal{X}_{\bm{\eta}'}).$$
For the same $S_{oa}=\{{\bm\eta}^1,\ldots,{\bm\eta}^m\}$ in generating $U_{oa}$, define 
$$\bar{V}=\frac{1}{m}\sum_{i=1}^m\bar{g}(\mathcal{X}_{\bm{\eta}^i}).$$

\begin{lemma}\label{lem:barfunctionproperty}
Some properties of $V$ and $\bar{V}$ are listed as follows.\\
{$(i)$} $\bar{V}$ is an unbiased estimator of ~$V$.\\
{$(ii)$} The bias of ~$V$ is of order $O(n^{-1})$ and ${\rm MSE}(V)={\rm MSE}(U_0)+O(n^{-2})$.\\
{$(iii)$} $U_{oa}$ is an unbiased estimator of ~$V$ and so also has bias $O(n^{-1})$.\end{lemma}

\noindent{\it Proof.} $(i)$ follows the unbiasedness of orthogonal arrays. $(ii)$ can be found in Proposition 3.5 in \cite{shao:2007} (page 211). $(iii)$ follows from \cite{owen:1992}.~~~$\square$
\begin{lemma}\label{lem:bar}
$$E(\bar{V}-V)^2\leq\frac{1}{m}\sum_{\bm u:|\bm u|>t}\left(\delta_{\bm u}^2+O(n^{-1})\right).$$
\end{lemma}
\noindent{\it Proof.}
Let $\delta^2_{\bm u}=\delta^2_{|\bm u|}$ and $\sigma^2_{\bm u}=\sigma^2_{|\bm u|}$.Change the $F$ in section 2.2 to $F_n$, we can define $dF_{n,{\bm u}}$, $g_{n,{\bm u}}$, $h_{n,{\bm u}}$, $\sigma^2_{n,{\bm u}}$ and $\delta^2_{n,{\bm u}}$ analogously and sequentially. Again, by substituting $\bar{g}$ for $g$, with $F_n$, we define $\bar{g}_{n,{\bm u}}$, $\bar{h}_{n,{\bm u}}$, $\bar{\sigma}^2_{n,{\bm u}}$ and $\bar{\delta}^2_{n,{\bm u}}$.
Adopt  (3.5) in Owen (1992) to $\bar{g}$, we have
$$E[(\bar{V}-V)^2|F_n]\leq\frac{1}{m}\sum_{\bm u:|\bm u|>t}\bar{\delta}^2_{n,\bm u}\leq\frac{1}{m}\sum_{\bm u:|\bm u|>t}{\delta}^2_{n,\bm u},$$
which leads to
$E(\bar{V}-V)^2=E(E[(\bar{V}-V)^2|F_n])\leq \frac{1}{m}\sum_{\bm u:|\bm u|>t}E{\delta}^2_{n,\bm u}$.
Consider $\sigma_{n,\bm u}^2=\int g^2_{n,\bm u}(x_1,\ldots,x_d)dF_n(x_1)\ldots dF_n(x_n)$, which can be further written as
$\int \left(\int g^2_{n,\bm u}dF_{n,\bm u^c}\right)^2dF_{n,\bm u}$.
This integer can be viewed as a V-statistic with the new kernel 
$g(x_1,\ldots,x_{|\bm u|},x_{|\bm u|+1},\ldots,x_d)\cdot$ $g(x_1,\ldots,x_{|\bm u|},x_{d+1},$ $\ldots,x_{2d-|\bm u|})$,
which estimates $\sigma_{\bm u}^2$ with bias $O(n^{-1})$.~~~$\square$


\begin{lemma}\label{lem:lusin} {\rm (Lusin's theorem)}

\noindent For any measurable function $g$ on $R^d$ and arbitrary $\epsilon>0$, there exists a continuous $g_{\epsilon}$ defined on $R^d$ with compact support such that $E|g-g_{\epsilon}|<\epsilon$.
\end{lemma}


\noindent{\large\bf Proof of Theorem \ref{thm:noassumptiononesample}.} Define $g_{F}(Z_1,\ldots,Z_d)=g(F^{-1}(Z_1),\ldots,F^{-1}(Z_d))$ such that $Z\sim U(0,1)$ and $F^{-1}(Z)\sim F$. With this new kernel $g_F$, the distribution of random variables $X$ is assumed to be the uniform distribution on $[0,1]$.

Write $U_{oa}-\Theta$ as $(U_{oa}-\bar{V})+(\bar{V}-V)+(V-\Theta)$.
Simple analysis reveals the following relationships among of $V_{oa}$, $\bar{V}$ and $V$.
Conditional on $F_n$, $V$ is constant and so $E(U_{oa}-\bar{V})(V-\Theta)=0$, $E(\bar{V}-V)(V-\Theta)=0$ since $E(U_{oa}-\bar{V})=E(\bar{V}-V)=0$. Conditional on both $V$ and $\bar{V}$, $E(U_{oa}-\bar{V})=0$ which indicates $E(U_{oa}-\bar{V})(\bar{V}-V)=0$.
Thus,
\begin{equation}\label{eq:simpledecomposition}
{\rm MSE}(U_{oa})=E(U_{oa}-\bar{V})^2+E(\bar{V}-V)^2+{\rm MSE}(V)
\end{equation}
whose last two terms have been addressed by Lemma \ref{lem:bar} and Lemma \ref{lem:barfunctionproperty}.
So we need to prove $E(U_{oa}-\bar{V})^2=o(m^{-1})$.
Since $U_{oa}$ and $\bar{V}$ always use the same $S_{oa}=\{\bm\eta^1,\ldots,\bm\eta^m\}$,
$$E(U_{oa}-\bar{V})^2=E\left(\frac{1}{m}\sum_{i=1}^mg(\mathcal{X}_{\bm{\eta}^i})-\bar{g}(\mathcal{X}_{\bm{\eta}^i})\right)^2.$$
For ${i_1}\neq {i_2}$ ($i_1,i_2\in \mathcal{Z}_m$), $E(g(\mathcal{X}_{\bm{\eta}^{i_1}})-\bar{g}(\mathcal{X}_{\bm{\eta}^{i_1}}))(g(\mathcal{X}_{\bm{\eta}^{i_2}})-\bar{g}(\mathcal{X}_{\bm{\eta}^{i_2}}))=0$. Denote $\bm\eta\sim\bm\eta'$ if $\bm\eta$ and $\bm\eta'$ belong to the same grid.
\begin{equation}\label{eq:funproof}
E(U_{oa}-\bar{V})^2\leq 2m^{-1}E[(g(\mathcal{X}_{\bm{\eta}})-{g}(\mathcal{X}_{\bm{\eta}'}))^2|{\bm \eta\sim\bm \eta'}].
\end{equation}

For any $M>0$, define
$g({\bm x},M)=\max\{\min\{g({\bm x}),M\},-M\}$.
Obviously, we have $\lim_{M\rightarrow\infty}g({\bm x},M)=g({\bm x})$, and dominated convergence theorem indicates
\begin{eqnarray}\label{eq:funproof1}
E[(g(\mathcal{X}_{\bm \eta})-g(\mathcal{X}_{\bm \eta'}))^2|{\bm \eta\sim\bm \eta'}]=\lim_{M\rightarrow\infty}E[(g(\mathcal{X}_{\bm \eta},M)-g(\mathcal{X}_{\bm \eta'},M))^2|{\bm \eta\sim\bm \eta'}].
\end{eqnarray}
Thus, for arbitrary $\epsilon>0$, we can find $M_{\epsilon}$ such that
\begin{eqnarray}\label{eq:funproof2}
&E[(g(\mathcal{X}_{\bm \eta})-g(\mathcal{X}_{\bm \eta'}))^2|{\bm \eta\sim\bm \eta'}]\leq E[(g(\mathcal{X}_{\bm \eta},M_{\epsilon})-g(\mathcal{X}_{\bm \eta'},M_{\epsilon}))^2|{\bm \eta\sim\bm \eta'}]+\epsilon.~~~~
\end{eqnarray}
Note that $\{X_1,\ldots,X_n\}$ are random, so is $\mathcal{X}_{\bm\eta}$.  Note that $Eg^2(X_1,\ldots,X_d)<\infty$.
We have $Eg^2(\mathcal{X}_{\bm\eta})<\infty$ and so $Eg(\mathcal{X}_{\bm\eta})<\infty$, which indicates $Eg^2(\mathcal{X}_{\bm \eta},M_{\epsilon})<\infty$ and $Eg(\mathcal{X}_{\bm\eta},M_{\epsilon})<\infty$.
From Lusin's theorem, there exists a continuous $g_{\epsilon,M_{\epsilon}}^*$ with compact support such that $E|g(\mathcal{X}_{\bm\eta},M)-g_{\epsilon,M_{\epsilon}}^*(\mathcal{X}_{\bm\eta})|<\epsilon M_{\epsilon}^{-1}$. Since $|g(\mathcal{X}_{\bm\eta},M_{\epsilon})|\leq M_{\epsilon}$,
\begin{eqnarray}\label{eq:funproof3}
&&E[(g(\mathcal{X}_{\bm \eta},M_{\epsilon})-g(\mathcal{X}_{\bm \eta'},M_{\epsilon}))^2|{\bm \eta\sim\bm \eta'}]\notag\\
&\leq&2M_{\epsilon}E[|g(\mathcal{X}_{\bm \eta},M_{\epsilon})-g(\mathcal{X}_{\bm \eta'},M_{\epsilon})||{\bm \eta\sim\bm \eta'}]\notag\\
&\leq&2M_{\epsilon}E|g(\mathcal{X}_{\bm\eta},M_{\epsilon})-g_{\epsilon,M_{\epsilon}}^*(\mathcal{X}_{\bm \eta})|+2M_{\epsilon}E|g(\mathcal{X}_{\bm \eta'},M_{\epsilon})-g_{\epsilon,M_{\epsilon}}^*(\mathcal{X}_{\bm\eta'})|+\notag\\
&&2M_{\epsilon}E[|g_{\epsilon,M_{\epsilon}}^*(\mathcal{X}_{\bm \eta})-g_{\epsilon,M_{\epsilon}}^*(\mathcal{X}_{\bm \eta'})||{\bm \eta\sim\bm \eta'}]\notag\\
&\leq&4\epsilon+2M_{\epsilon}E[|g_{\epsilon,M_{\epsilon}}^*(\mathcal{X}_{\bm \eta})-g_{\epsilon,M_{\epsilon}}^*(\mathcal{X}_{\bm \eta'})||{\bm \eta\sim\bm \eta'}]
\end{eqnarray}
Note that $g_{\epsilon,M_{\epsilon}}^*$ has compact support and so is uniformly continuous. There exists $\Delta(M_{\epsilon}^{-1}\epsilon)$ such that $|g_{\epsilon,M_{\epsilon}}^*(\mathcal{X}_{\bm \eta})-g_{\epsilon,M_{\epsilon}}^*(\mathcal{X}_{\bm \eta'})|\leq \epsilon M_{\epsilon}^{-1}$ as long as $||\mathcal{X}_{\bm \eta}-\mathcal{X}_{\bm \eta'}||_2\leq \Delta(M_{\epsilon}^{-1}\epsilon)$.
Define $$\mathcal{A}=\{|\mathcal{X}_{\eta_{j}}-\mathcal{X}_{\eta'_{j}}|\geq d^{-1}\Delta(M_{\epsilon}^{-1}\epsilon)~{\rm for~some}~j\in\mathcal{Z}_d\},$$ 
with $P(\mathcal{A})\leq\sum_{j=1}^{d}P\{|\mathcal{X}_{ \eta_{j}}-\mathcal{X}_{\eta'_{j}}|\geq d^{-1}\Delta(M_{\epsilon}^{-1}\epsilon)\}$, and $||\mathcal{X}_{\bm \eta}-\mathcal{X}_{\bm \eta'}||_2\leq \Delta(M_{\epsilon}^{-1}\epsilon)$ on $\mathcal{A}^c$.
\begin{eqnarray}
&&2M_{\epsilon}E[|g_{\epsilon,M_{\epsilon}}^*(\mathcal{X}_{\bm \eta})-g_{\epsilon,M_{\epsilon}}^*(\mathcal{X}_{\bm \eta'})||{\bm \eta\sim\bm \eta'}]\notag\\
&=&2M_{\epsilon}P(\mathcal{A}^c)E[|g_{\epsilon,M_{\epsilon}}^*(\mathcal{X}_{\bm \eta})-g_{\epsilon,M_{\epsilon}}^*(\mathcal{X}_{\bm \eta'})||{\bm\eta\sim\bm \eta'},\mathcal{A}^c]\notag\\
&&+2M_{\epsilon}P(\mathcal{A})E[|g_{\epsilon,M_{\epsilon}}^*(\mathcal{X}_{\bm \eta})-g_{\epsilon,M_{\epsilon}}^*(\mathcal{X}_{\bm \eta'})||{\bm \eta\sim\bm \eta'},\mathcal{A}]\notag\\
&\leq&2\epsilon+4M^2_{\epsilon}\sum_{k=1}^{d}P\{|\mathcal{X}_{ \eta_{j}}-\mathcal{X}_{\eta'_{j}}|\geq d^{-1}\Delta(M_{\epsilon}^{-1}\epsilon)\}
\end{eqnarray}
Now we give the relationship among several events. For $j\in\mathcal{Z}_d$ and $\bm\eta\sim\bm\eta'$,
\begin{eqnarray}
&&\{|\mathcal{X}_{\eta_{j}}-\mathcal{X}_{\eta'_{j}}|\geq d^{-1}\Delta(M_{\epsilon}^{-1}\epsilon)\}\notag\\
&=&\{|\mathcal{X}_{\eta_{j}}-F_{n}(\mathcal{X}_{\eta_{j}})+F_{n}(\mathcal{X}_{\eta_{j}})-F_{n}(\mathcal{X}_{\eta'_{j}})+F_{n}(\mathcal{X}_{\eta'_{j}})-\mathcal{X}_{\eta'_{j}}|\geq d^{-1}\Delta(M_{\epsilon}^{-1}\epsilon)\}\notag\\
&\subseteq&\{\sup_{x\in(0,1)}|x-F_{n}(x)|\geq \frac{1}{3d}\Delta(M_{\epsilon}^{-1}\epsilon)\}\cup\{F_{n}(\mathcal{X}_{\eta_{j}})-F_{n}(\mathcal{X}_{\eta'_{j}})\geq \frac{1}{3d}\Delta(M_{\epsilon}^{-1}\epsilon)\}\notag
\end{eqnarray}
Note that $\bm\eta\sim\bm\eta'$, as $L\rightarrow\infty$, $P(\{F_{n}(\mathcal{X}_{\eta_{j}})-F_{n}(\mathcal{X}_{\eta'_{j}})\geq \frac{1}{3d}\Delta(M_{\epsilon}^{-1}\epsilon)\})\rightarrow 0$.
Dvoretzky-Kiefer-Wolfowitz inequality reveals $P\left(\sup_{x\in (0,1)}|F_{n}(x)-x|\geq \epsilon\right)\leq \exp(-2n\epsilon^2)$. So we immediately have
$P(\{|\mathcal{X}_{\eta_{j}}-\mathcal{X}_{\eta'_{j}}|\geq d^{-1}\Delta(M_{\epsilon}^{-1}\epsilon)\})\rightarrow 0$
as $n, L\rightarrow\infty$, and we can find $n_{\epsilon}$ and $L_{\epsilon}$ such that
\begin{equation}\label{eq:funproof4}
P(\{|\mathcal{X}_{\eta_{j}}-\mathcal{X}_{\eta'_{j}}|\geq d^{-1}\Delta(M_{\epsilon}^{-1}\epsilon)\})\leq (4dM^2_{\epsilon})^{-1}\epsilon
\end{equation}
as long as $n\geq n_{\epsilon}$ and $L\geq L_{\epsilon}$.

Finally, by combining (\ref{eq:funproof1})-(\ref{eq:funproof4}), we know that for arbitrary $\epsilon>0$, we can find $n_{\epsilon}$ and $L_{\epsilon}$ such that
$E[(g(\mathcal{X}_{\bm\eta})-g(\mathcal{X}_{\bm \eta'}))^2|{\bm \eta\sim\bm \eta'}]\leq 8\epsilon$,
as long as $n\geq n_{\epsilon}$ and $L\geq L_{\epsilon}$.
That means
\begin{equation}\label{eq:funprooffinal}
E[(g(\mathcal{X}_{\bm \eta})-g(\mathcal{X}_{\bm \eta'}))^2|{\bm \eta\sim\bm \eta'}]\rightarrow 0
\end{equation}
as $n,L\rightarrow\infty$.
Theorem \ref{thm:noassumptiononesample} is concluded by submitting (\ref{eq:funprooffinal}) into (\ref{eq:funproof}) and combining (\ref{eq:funproof}) with (\ref{eq:simpledecomposition}), Lemma \ref{lem:barfunctionproperty}($ii$) and Lemma \ref{lem:bar}.~~~$\square$
\vspace{.5cm}

\vspace{.5cm}

\noindent{\large\bf Proof of Theorem \ref{thm:lipschitzonesample}.} 
There exists $c>0$ such that density function $f(\cdot)>c$ on $[a,b]$, and $|F(x_1)-F(x_2)|\geq c|x_1-x_2|$ for $x_1,x_2\in [a,b]$.
In  (\ref{eq:simpledecomposition}), we only analyze $E(U_{oa}-\bar{V})^2$ since the rest two terms are given by Lemma \ref{lem:barfunctionproperty}$(ii)$ and Lemma \ref{lem:bar}.
Dvoretzky-Kiefer-Wolfowitz inequality reveals $P\left(\sup_{x\in R}|F_{n}(x)-F(x)|\geq \epsilon\right)\leq \exp(-2n\epsilon^2)$.
By taking $\epsilon=[\log(n)n^{-1}]^{1/2}$, we have
\begin{eqnarray}\label{eqDKW}
P\left(\mathcal{A}\right)\leq\exp(-2\log n)=O(n^{-2}),\notag
\end{eqnarray}
where $\mathcal{A}=\{\sup_{x\in R}|F_{n}(x)-F(x)|\geq n^{-1/2}\log^{1/2}(n)\}$.
Since $g$ is continuous and $F$ is bounded, we can find $M>0$ such that $|g|\leq M$ and so $|U_{oa}|,|\bar{V}|\leq M$.
\begin{eqnarray}
&&E[(g(\mathcal{X}_{\bm{\eta}})-{g}(\mathcal{X}_{\bm{\eta}'}))^2|{\bm \eta\sim\bm \eta'}]\notag\\
&=&P(\mathcal{A})E[(g(\mathcal{X}_{\bm{\eta}})-{g}(\mathcal{X}_{\bm{\eta}'}))^2|{\bm \eta\sim\bm \eta'},\mathcal{A}]+P(\mathcal{A}^c)E[(g(\mathcal{X}_{\bm{\eta}})-{g}(\mathcal{X}_{\bm{\eta}'}))^2|{\bm \eta\sim\bm \eta'},\mathcal{A}^c]\notag\\
&\leq& M^2n^{-2}+E[(g(\mathcal{X}_{\bm{\eta}})-{g}(\mathcal{X}_{\bm{\eta}'}))^2|{\bm \eta\sim\bm \eta'},\mathcal{A}^c]\notag
\end{eqnarray}
The analysis of $E[(U_{oa}-\bar{V})^2|\mathcal{A}^c]$ is as follows. On $\mathcal{A}^c$, we have, for $1\leq k_1,k_2\leq nL^{-1}$,
\begin{eqnarray}
&&c|X_{((l-1)nL^{-1}+k_1)}-X_{((l-1)nL^{-1}+k_2)}|\notag\\
&\leq& |F(X_{((l-1)nL^{-1}+k_1)})-F(X_{((l-1)nL^{-1}+k_2)})|\notag\\
&\leq& |F_{n}(X_{((l-1)nL^{-1}+k_1)})-F_{n}(X_{((l-1)nL^{-1}+k_2)})|+2n^{-1/2}\log^{1/2}n\notag\\
&\leq& L^{-1}+2n^{-1/2}\log^{1/2}n.\notag
\end{eqnarray}
Since $g$ is Lipschitz continuous, we know
$(g(\mathcal{X}_{\bm\eta})-g(\mathcal{X}_{\bm\eta'}))^2=O(L^{-2}+n^{-1}\log n)$
for any $\bm\eta\sim\bm\eta'$.
Then we have $E[(g(\mathcal{X}_{\bm\eta})-g(\mathcal{X}_{\bm\eta'}))^2|{\bm\eta\sim\bm\eta'},\mathcal{A}^c]=O(L^{-2}+n^{-1}\log n)$. With this equation, Theorem \ref{thm:boundedvaronesample} is the direct result of (\ref{eq:simpledecomposition}) (\ref{eq:funproof}), Lemma \ref{lem:barfunctionproperty}$(ii)$, Lemma \ref{lem:bar}.~~~~$\square$

\vspace{.5cm}

\vspace{.5cm}

\noindent{\large\bf Proof of Theorem \ref{thm:tradeoff}.}
For convenience, we simply write $g_F$ as $g$ in this proof.
In (\ref{eq:simpledecomposition}), we only analyze $E(U_{oa}-\bar{V})^2$ since the rest two terms are given by Lemma \ref{lem:barfunctionproperty} $(ii)$ and Lemma \ref{lem:bar}.
Each row of the matrix $A$ generated in step 1 follows the uniform distribution on $\mathcal{Z}_L^d$ since the permutation in each column of $A_0$ is independent. Thus,
$$E(U_{oa}-\bar{V})^2=E\left(\frac{1}{m}\sum_{i=1}^mg(\mathcal{X}_{\bm{\eta}^i})-\bar{g}(\mathcal{X}_{\bm{\eta}^i})\right)^2=\frac{1}{mL^d}\sum_{\bm a\in \mathcal{Z}_L^d}E[(g(\mathcal{X}_{\bm{\eta}})-\bar{g}(\mathcal{X}_{\bm{\eta}}))^2|\bm\eta\in\mathcal{G}_{\bm a}].$$
Analysis is now focused on $E[(g(\mathcal{X}_{\bm{\eta}})-\bar{g}(\mathcal{X}_{\bm{\eta}}))^2|\bm\eta\in\mathcal{G}_{\bm a}]$ for every $\bm a\in\mathcal{Z}_L^d$. Let $X_{(0)}=0$ and $X_{(n+1)}=1$. For $l\in\mathcal{Z}_L$, given $X_{((l-1)nL^{-1})}$ and $X_{(lnL^{-1}+1)}$, $X_{((l-1)nL^{-1}+1)},\ldots,X_{(lnL^{-1})}$ has the same distribution as the order statistic of $L$ samples following the uniform distribution on $[X_{((l-1)nL^{-1})},X_{(lnL^{-1}+1)}]$.
For $\mathcal{A}=\{\sup_{x\in R}|F_{n}(x)-F(x)|\geq n^{-\frac{1-c}{2}}\}$, Dvoretzky-Kiefer-Wolfowitz inequality reveals $P(\mathcal{A})=\exp(-2n^c)$.
On $\mathcal{A}^c$, we have $(X_{(lnL^{-1}+1)}-X_{((l-1)nL^{-1})})/L\rightarrow 1$ as $n\rightarrow\infty$. The analysis is now focused on
$E[(g(\mathcal{X}_{\bm{\eta}})-\bar{g}(\mathcal{X}_{\bm{\eta}}))^2|\bm\eta\in\mathcal{G}_{\bm a},\mathcal{A}^c]$.
For this given $\bm a$, define $\mathcal{X}_0=(X_{0,1},\ldots,X_{0,d})$ where $X_{0,j}=\frac{L}{n}\sum_{\eta\in G_{a_j}}X_{\eta}$ and so $\sum_{\eta\in G_{a_j}}(X_{\eta}-X_{0,j})=0$. Adopt the Taylor expansion on $\mathcal{X}_0$, we have 
\begin{equation}\notag
g(\mathcal{X}_{\bm\eta})=g(\mathcal{X}_0)+\sum_{j=1}^d\frac{\partial g}{\partial x_j}\bigg|_{X_{0,j}}(X_{\eta_j}-X_{0,j})+O(L^{-2})~~{\rm and}~~\bar{g}(\mathcal{X}_{\bm\eta})=g(\mathcal{X}_0)+O(L^{-2}).
\end{equation}
\begin{eqnarray}
&&E[(g(\mathcal{X}_{\bm{\eta}})-\bar{g}(\mathcal{X}_{\bm{\eta}}))^2|\bm\eta\in\mathcal{G}_{\bm a},\mathcal{A}^c]\notag\\
&=&E\left[\left(\sum_{j=1}^d\frac{\partial g}{\partial x_j}\bigg|_{X_{0,j}}\cdot(X_{\eta_j}-X_{0,j})+O(L^{-2})\right)^2|\bm\eta\in\mathcal{G}_{\bm a},\mathcal{A}^c\right]\notag\\
&=&o(L^{-2})+\sum_{j=1}^dE\left[\left(\frac{\partial g}{\partial x_j}\bigg|_{X_{0,j}}\cdot(X_{\eta_j}-X_{0,j})\right)^2|\bm\eta\in\mathcal{G}_{\bm a},\mathcal{A}^c\right]\notag\\
&=&o(L^{-2})+\sum_{j=1}^d\left(\frac{\partial g}{\partial x_j}\bigg|_{X_{0,j}}\right)^2\frac{1}{12L^2}.\notag
\end{eqnarray}
And then we have
\begin{eqnarray}
&&E(U_{oa}-\bar{V})^2=\frac{1}{mL^d}\sum_{\bm a\in \mathcal{Z}_L^d}E[(g(\mathcal{X}_{\bm{\eta}})-\bar{g}(\mathcal{X}_{\bm{\eta}}))^2|\bm\eta\in\mathcal{G}_{\bm a}]\notag\\
&=&\frac{1}{12mL^2}\sum_{j=1}^d\left(\frac{1}{L^d}\sum_{\bm a\in \mathcal{Z}_L^d}\left(\frac{\partial g}{\partial x_j}\bigg|_{X_{0,j}}\right)^2\right)+o\left(\frac{1}{mL^2}\right)\notag\\
&=&\frac{1}{12mL^2}\sum_{j=1}^dE\left(\frac{\partial g}{\partial x_j}\right)^2+o\left(\frac{1}{mL^2}\right),\notag
\end{eqnarray}
Then Theorem \ref{thm:tradeoff} is the direct result of (\ref{eq:simpledecomposition}), Lemma \ref{lem:barfunctionproperty}($ii$) and Lemma \ref{lem:bar}.~~~~$\square$

\vspace{.5cm}

\noindent{\large\bf Proof of Theorem \ref{thm:debias}.}
Consider the $m$ rows of $A$, $\bm a_1,\ldots,\bm a_m$, generated in the step 1 of the construction in section 2.1. For any $\bm a\in\mathcal{Z}_L^d$, the random permutation in generating $\bm a_1,\ldots,\bm a_m$ reveals that $P(\bm a_1=\bm a)=L^{-d}$. Given $F_n$,
\begin{eqnarray}
E(\tilde{U}_{oa}|F_n)&=&E\frac{1}{m}\sum_{i=1}^m\omega_{\bm{\eta}^i}g(\mathcal{X}_{\bm \eta^{i}})=E\omega_{\bm{\eta}^1}g(\mathcal{X}_{\bm \eta^{1}})\notag\\
&=&\sum_{\bm a\in \mathcal{Z}_L^d}L^{-d}E_{\bm\eta\in\mathcal{G}_{\bm a}}\omega_{\bm{\eta}}g(\mathcal{X}_{\bm \eta})=\sum_{\bm a\in \mathcal{Z}_L^d}\frac{|\mathcal{G}_{{\bm a}^i}\cap S_0^*|}{|S_0^*|}E_{\bm\eta\in\mathcal{G}_{\bm a}}g(\mathcal{X}_{\bm\eta})\notag\\
&=&\sum_{\bm a\in \mathcal{Z}_L^d}\frac{|\mathcal{G}_{{\bm a}^i}\cap S_0^*|}{|S_0^*|}\left(\frac{1}{|\mathcal{G}_{{\bm a}^i}\cap S_0^*|}\sum_{\bm\eta\in\mathcal{G}_{\bm a}}g(\mathcal{X}_{\bm\eta})\right)=\frac{1}{|S_0^*|}\sum_{\bm\eta\in S_0^*}g(\mathcal{X}_{\bm\eta})=U_0.\notag
\end{eqnarray}
Since $U_0$ is unbiased, so is $\tilde{U}_{oa}$. This proves the unbiasedness of $\tilde{U}_{oa}$. The MSE of $\tilde{U}_{oa}$ can be similar analyzed as Theorem \ref{thm:noassumptiononesample}, and so is omitted here.~~~~$\square$


\bibhang=1.7pc
\bibsep=2pt
\fontsize{9}{14pt plus.8pt minus .6pt}\selectfont
\renewcommand\bibname{\large \bf References}
\expandafter\ifx\csname
natexlab\endcsname\relax\def\natexlab#1{#1}\fi
\expandafter\ifx\csname url\endcsname\relax
  \def\url#1{\texttt{#1}}\fi
\expandafter\ifx\csname urlprefix\endcsname\relax\def\urlprefix{URL}\fi

\lhead[\footnotesize\thepage\fancyplain{}\leftmark]{}\rhead[]{\fancyplain{}\rightmark\footnotesize{} }

\vskip .65cm
\noindent
School of Mathematics and Statistics, Beijing Institute of Technology
\vskip 2pt
\noindent
E-mail: kongsunday@163.com
\vskip 2pt
\noindent
Department of Business Analytics and Statistics, University of Tennessee
\vskip 2pt
\noindent
E-mail: wzheng9@utk.edu
\newpage

\fontsize{12}{14pt plus.8pt minus .6pt}\selectfont \vspace{0.8pc}
\markboth{\hfill{\footnotesize\rm XIANGSHUN KONG AND WEI ZHENG} \hfill}
{\hfill {\footnotesize\rm SUPPLEMENTARY MATERIAL} \hfill}
\centerline{\large\bf DESIGN BASED INCOMPLETE U-STATISTICS}
\vspace{.4cm} \centerline{Xiangshun Kong$^1$, Wei Zheng$^2$} \vspace{.4cm} \centerline{\it
    $^1$Beijing Institute of Technology and $^2$University of Tennessee} \vspace{.55cm} \fontsize{9}{11.5pt plus.8pt minus
.6pt}\selectfont
\fontsize{12}{14pt plus.8pt minus .6pt}\selectfont \vspace{0.8pc}

\noindent{\large\bf{Generalization of Theorem \ref{thm:lipschitzonesample}.}}

The following conditions on $g$ or $F$ will be needed by Theorem \ref{thm:lipschitzonesample} in Section 2 and Theorems \ref{thm:clipschitzonesample}--\ref{thm:linearcombination} in this section.
\begin{itemize}
\item[($g.1$)] {\it Lipschitz continuous}: The function, $g:R^d\rightarrow R$, is said to be Lipschitz continuous if there exists a constant $c>0$ such that $|g({\bf a_1})-g({\bf a_2})|\leq c||{\bf a_1}-{\bf a_2}||_2$ for any ${\bf a_1},{\bf a_2}\in R^d$. Example: First-order polynomial functions.
\item[($g.2$)] {\it Order-$p$ continuous}: The function, $g:R^d\rightarrow R$, is said to be order-$p$ continuous if there exists a constant $c>0$ and $\phi_p({\bf a_1}-{\bf a_2})\leq c+\max^{p}(||{\bf a_1}||_2,||{\bf a_2}||_2)$ for any ${\bf a_1},{\bf a_2}\in R^d$ such that $|g({\bf a_1})-g({\bf a_2})|\leq\phi({\bf a_1},{\bf a_2})||{\bf a_1}-{\bf a_2}||_2$ for any ${\bf a_1},{\bf a_2}\in R^d$. \\
Example: All polynomial functions.
\item[($g.3$)] {\it Uniformly bounded-variation}: For a real valued function $f:R\rightarrow R$, the total variation of $f$ is defined as
$V_R(f)=\sup_{p>0}\sup_{-\infty< c_1,\ldots,c_p<\infty}\sum_{i=1}^{p-1}|f(c_{i+1})-f(c_i)|.$ The function, $g:R^d\rightarrow R$, is said to be uniformly bounded-variation if there exists a constant $c>0$ such that $V_R(g(\cdot,x_2,\ldots,x_d))<c$ for any $(x_2,\ldots,x_d)\in R^{d-1}$. \\
Example: Linear combinations of sign functions, e.g. $g(x_1,x_2)={\rm sign}(x_1x_2)+{\rm sign}(x_1+x_2)$.
\item[($F$)] {\it Light-tailed} distribution: The distribution of a random variable $X$ is said to be light-tailed if there exists constants $c,c_1>0$ such that $P(|X|>x)\leq e^{-cx}$ for all $x>c_1$. Example: Normal distribution, exponential distribution, and truncated distributions.
\end{itemize}

\begin{lemma}\label{lem:moment}
Suppose $F$ is light-tailed. Let $X_{\max}=\max\{|X_1|,\ldots,|X_n|\}$. Then, for arbitrary $a>0$ with $n\rightarrow\infty$, we have
$$EX_{\max}^a=O(\log n)^a.$$
\end{lemma}
\noindent{\it Proof.} Since the distribution is light-tailed, we have $P(|X|>x)\leq e^{-cx}$ for any $|x|>c_0$, where $c$ and $c_0$ are two fixed positive numbers.
\begin{eqnarray}
E(X_{\max})^a&=&\int_{x>0}ax^{a-1}P(X_{\max}>x)dx\notag\\
&\leq&\int_0^{2c^{-1}\log n}ax^{a-1}dx+\int_{2c^{-1}\log n}^{\infty}ax^{a-1}P(X_{\max}>x)dx\notag\\
&=& O(\log n)^a+\int_{2c^{-1}\log n}^{\infty}ax^{a-1}P(X_{\max}>x)dx\notag\\
&=& O(\log n)^a+\int_{2c^{-1}\log n}^{\infty}ax^{a-1}ne^{-cx}dx~=~O(\log n)^a+O(1).\notag~~~\square
\end{eqnarray}

\begin{lemma}\label{lem:voa}
Suppose $(i)$ $g$ is order-$p$ continuous, and $(ii)$ $F$ is light-tailed. We have
$$E(U_{oa}-\bar{V})^2= O\left(\frac{1}{mL}(\log n)^{2p+2}\right).$$
\end{lemma}
\noindent{\it Proof.} Let $X_{\max}=\max\{|X_1|,\ldots,|X_n|\}$. For $l\in\mathcal{Z}_L$, define
$d_{l}=\max\{|X_{i_1}-X_{i_2}|:i_1,i_2\in G_l\}.$
Since $g$ is order-$p$ continuous, for $\bm{\eta}\sim\bm{\eta}'$ in $\mathcal{G}_{\bm a}$,
$|g(\mathcal{X}_{\bm{\eta}})-g(\mathcal{X}_{\bm{\eta}'})|\leq (c_1+X_{\max}^{p})d^{1/2}d_l$,
and so
$|g(\mathcal{X}_{\bm{\eta}})-g(\mathcal{X}_{\bm{\eta}'})|^2\leq (c_1+X_{\max}^{p})^2\cdot d\cdot\sum_{j=1}^{d}d^2_{a_{j}}$.

Since $U_{oa}$ and $\bar{V}$ always use the same $S_{oa}=\{\bm\eta^1,\ldots,\bm\eta^m\}$, we have
$$E(U_{oa}-\bar{V})^2=E\left(\frac{1}{m}\sum_{i=1}^m(g(\mathcal{X}_{\bm{\eta}^i})-\bar{g}(\mathcal{X}_{\bm{\eta}^i}))\right)^2.$$
For $\bm{i}_1\neq \bm{i}_2$, $E(g(X_{\bm{\eta}^{i_1}})-\bar{g}(X_{\bm{\eta}^{i_1}}))(g(X_{\bm{\eta}^{i_2}})-\bar{g}(X_{\bm{\eta}^{i_2}}))=0$.
\begin{eqnarray}\label{eq:innerdifference}
E(U_{oa}-\bar{V})^2&=&m^{-2}E\sum_{i=1}^m(g(\mathcal{X}_{\bm{\eta}^i})-\bar{g}(\mathcal{X}_{\bm{\eta}^i}))^2\notag\\
&\leq& m^{-2}E\sum_{i=1}^{m}(c_1+X_{\max}^{p})^2\cdot d\cdot\sum_{j=1}^{d}d^2_{a^i_{j}}\notag
\end{eqnarray}
Since $\sum_{l=1}^Ld_{l}\leq 2X_{\max}$, we have $\sum_{l=1}^Ld_{l}^2\leq 4X_{\max}^2$. Using Lemma \ref{lem:moment}, we have
\begin{eqnarray}
E(U_{oa}-\bar{V})^2&\leq& m^{-2}dE\left((c_1+X_{\max}^{p})^2\sum_{i=1}^{m}\sum_{j=1}^{d}d^2_{a^i_{j}}\right)=m^{-2}dE\left((c_1+X_{\max}^{p})^2\sum_{j=1}^{d}\sum_{i=1}^{m}d^2_{a^i_{j}}\right)\notag\\
&=& m^{-2}dE\left((c_1+X_{\max}^{p})^2\sum_{j=1}^{d}mL^{-1}4X_{\max}^2\right)=O\left(\frac{1}{mL}(\log n)^{2p+2}\right).~~~\square\notag
\end{eqnarray}

\begin{theorem}\label{thm:clipschitzonesample}
Suppose $(i)$ The kernel function $g$ is order-$p$ continuous, and $(ii)$ $F$ is light-tailed. For $U_{oa}$ based on $OA(m,d,L,t)$, we have
\begin{eqnarray}\label{eq:clipschitzonesample}
{\rm MSE}(U_{oa})={\rm MSE}(U_0)+\frac{R(t)}{m}+O\left(\frac{(\log n)^{2p+2}}{mL}\right)+O\left(\frac{1}{n^2}\right).
\end{eqnarray}
\end{theorem}
\noindent{\it Proof}. This is the direct result of (\ref{eq:simpledecomposition}), Lemma \ref{lem:barfunctionproperty}($ii$), Lemmas \ref{lem:bar} and \ref{lem:voa}.~~~~$\square$

\begin{theorem}\label{thm:boundedvaronesample}
Suppose the kernel function $g$ has uniformly bounded variation. For $U_{oa}$ based on $OA(m,d,L,t)$, we have
\begin{eqnarray}\label{eq:boundedvaronesample}
{\rm MSE}(U_{oa})={\rm MSE}(U_0)+\frac{R(t)}{m}+O\left(\frac{1}{mL}\right)+O\left(\frac{1}{n^2}\right).
\end{eqnarray}
\end{theorem}

\noindent{\it Proof.}
From (\ref{eq:simpledecomposition}), Lemma \ref{lem:barfunctionproperty}($ii$) and Lemma \ref{lem:bar}, we only need to prove $E(U_{oa}-\bar{V})^2=O(m^{-1}L^{-1})$. First, we introduce some notations that will be used only in the proof of this theorem. Given the order statistic of $\{X_1,\ldots,X_n\}$ denoted by $X_{(1)},\ldots,X_{(n)}$, for $l=1,\ldots,L$ and $(x_2,\ldots,x_d)\in R^{d-1}$, define
$D(l|x_2,\ldots,x_k)=\max_{(l-1)nL^{-1}<i_1<i_2\leq l\cdot nL^{-1}}$ $|g(X_{(i_1)},$ $x_2,\ldots,x_k)-g(X_{(i_2)},x_2,\ldots,x_k)|$.
Since $g$ has uniformly bounded variation, $g$ is bounded, say $|g|\leq M$.
\begin{eqnarray}
E[(g(\mathcal{X}_{\bm{\eta}})-{g}(\mathcal{X}_{\bm{\eta}'}))^2|{\bm \eta\sim\bm \eta'}]&=&L^{-d}\sum_{{\bm a}\in\mathcal{Z}_L^d}|\mathcal{G}_{\bm a}|^{-2}\sum_{{\bm\eta}\in \mathcal{G}_{\bm a}}\sum_{{\bm\eta'}\in \mathcal{G}_{\bm a}}(g(\mathcal{X}_{\bm \eta})-{g}(\mathcal{X}_{\bm \eta'}))^2\notag\\
&\leq&2ML^{-d}|\mathcal{G}_{\bm a}|^{-2}\sum_{{\bm a}\in\mathcal{Z}_L^d}\sum_{{\bm\eta}\in \mathcal{G}_{\bm a}}\sum_{{\bm\eta'}\in \mathcal{G}_{\bm a}}|g(\mathcal{X}_{\bm \eta})-{g}(\mathcal{X}_{\bm \eta'})|.\notag
\end{eqnarray}
Note that $g(\mathcal{X}_{\bm \eta})-{g}(\mathcal{X}_{\bm \eta'})$ can be written as the summation of the difference in changing each element of $\mathcal{X}_{\bm\eta}=(X_{\eta_1},\ldots,X_{\eta_d})$ to $\mathcal{X}_{\bm\eta'}=(X_{\eta'_1},\ldots,X_{\eta'_d})$ one by one as follows.
\begin{eqnarray}
&&|g(\mathcal{X}_{\bm \eta})-{g}(\mathcal{X}_{\bm \eta'})|\notag\\
&=&|g(X_{\eta_{1}},X_{\eta_{2}},\cdots)-g(X_{\eta'_{1}},X_{\eta_{2}},\cdots)|+|g(X_{\eta'_{1}},X_{\eta_{2}},X_{\eta_{3}},\cdots)-g(X_{\eta'_{1}},X_{\eta'_{2}},X_{\eta_{3}},\cdots)|\notag\\
&+&\cdots+|g(X_{\eta'_{1}},X_{\eta'_{2}},X_{\eta'_{3}},\cdots,X_{\eta'_{d-1}},X_{\eta_d})-g(X_{\eta'_{1}},X_{\eta'_{2}},X_{\eta'_{3}},\cdots,X_{\eta'_{d-1}},X_{\eta'_d})|\notag\\
&\leq& D(a_1 | X_{\eta_2},\ldots,X_{\eta_d})\notag+D(a_2 | X_{\eta'_1},X_{\eta_3},\ldots,X_{\eta_d})+\cdots+D(a_d | X_{\eta'_1},X_{\eta'_3},\ldots,X_{\eta'_{d-1}})
\end{eqnarray}

For orthogonal arrays, we can separate $\sum_{{\bm a}\in\mathcal{Z}_L^d}\sum_{{\bm\eta}\in \mathcal{G}_{\bm a}}\sum_{{\bm\eta'}\in \mathcal{G}_{\bm a}}D(a_1 | X_{\eta_2},\ldots,X_{\eta_d})$ into $|\mathcal{Z}_L^d||\mathcal{G}_{\bm a}|^{2}/L$ groups such that each group contains $L$ elements whose summation is control by the total variation $c>0$. So we have $$\sum_{{\bm a}\in\mathcal{Z}_L^d}\sum_{{\bm\eta}\in \mathcal{G}_{\bm a}}\sum_{{\bm\eta'}\in \mathcal{G}_{\bm a}}D(a_1 | X_{\eta_2},\ldots,X_{\eta_d})\leq cL^d|\mathcal{G}_{\bm a}|^2/L.$$ Similarly analyzing the $D(a_2 | X_{\eta'_1},X_{\eta_3},\ldots,X_{\eta_d})$, $\ldots$, $D(a_d | X_{\eta'_1},X_{\eta'_3},\ldots,X_{\eta'_{d-1}}$, we have
$E[(g(\mathcal{X}_{\bm{\eta}})-{g}(\mathcal{X}_{\bm{\eta}'}))^2|{\bm \eta\sim\bm \eta'}]=O(L^{-1})$ and so $E(U_{oa}-\bar{V})^2=O(m^{-1}L^{-1})$. Theorem \ref{thm:boundedvaronesample} is the direct result of (\ref{eq:simpledecomposition}), Lemma \ref{lem:barfunctionproperty}$(ii)$, Lemma \ref{lem:bar}. ~~~~$\square$

\begin{theorem}\label{thm:linearcombination}
Suppose $(i)$ The kernel function $g$ is a linear combination of some order-$p$ continuous functions and some uniformly bounded-variation functions, and $(ii)$ $F$ is light-tailed. Then (\ref{eq:clipschitzonesample}) still holds with $L^2\leq n(\log n)^{-1}$.
\end{theorem}

\noindent{\it Proof}. This is the direct result of Theorems \ref{thm:clipschitzonesample} and \ref{thm:boundedvaronesample}.

\noindent{\large\bf{Choosing $L$ and $t$.}}

From Eq(2.13) of Theorem 3 in the manuscript and the relation $m=\lambda L^t$, we know that the trade-off between $L$ and $t$ depends on the variance of each component in the Heoffding's decomposition, i.e., $\delta^2_j$, $j=1,\ldots,d$. We shall give these variances a estimator $\hat{\delta}^2_j$.
Using Eq(2.13) with $R(t)$ and ${E}\gamma^2(X_1,\ldots,X_d)$ being estimated as a function of $\hat{\delta}^2_j$, we should choose the combination of $L$ and $t$ which minimizes
$$\phi(L,t)=\frac{\hat{R}(t)}{m}+\frac{d}{12mL^2}\hat{E}\gamma^2(X_1,\ldots,X_d),$$
where $\hat{R}(t)$ and $\hat{E}\gamma^2(X_1,\ldots,X_d)$ are functions of $\hat{\delta}^2_j$'s.

Now we provide two methods for generating $\hat{\delta}^2_j$. (1) When the Heoffding's decomposition is easy to calculate, one can write down the analytical expression and give a direct estimation of $\delta^2_j$'s. (2) We can use a bootstrap approach for $\hat{\delta}^2_j$'s. With a small sample size $n'\ll n$, it is easy to bootstrap ${\rm MSE}(U_0)$ (the complete U-statistic). For details of the bootstrap approach, we may refer to Marie Huskova and Paul Janssen (1993a,b). Now, let us review the formula of ${\rm MSE}(U_0)$:
$${\rm MSE}(U_0)=\binom{n}{d}^{-1}\sum_{j=1}^{d}\binom{d}{j}\binom{n-d}{d-j}\sigma_{j}^2~=~\sum_{j=1}^d\binom{d}{j}^2\binom{n}{j}^{-1}\delta_j^2.$$
Usually, with at most $d$ different $n'(>d)$, we can generate linear equations of $\delta_j^2$ based on the $d$ different ${\widehat{\rm MSE}}(U_0)$ based on the bootstrap approach. And the solution of these linear equations can be used as the estimation of $\hat{\delta}^2_j$'s.

For the second method, we now use the setup in Example \ref{example:manycomparison} for illustration. For convenience, we set $n=10^4$ and $m=10^6$. The two choices of the combination of $L$ and $t$ is $(L=100,t=3)$ and $(L=1000,t=2)$. We use bootstrap method to estimate the variance of the complete U-statistic with $n'=4,5,6$. The subsample size $n'$ is so small that the computational burden of the bootstrapped complete U-statistic, i.e., $\binom{n'}{3}$ is negligible. Simulation reveals that $\hat{\delta}_1=0.0557$, $\hat{\delta}_2=0.00217$ and $\hat{\delta}_3=1.06257$. Simple analysis reveals that $t=3$ shall work better than $t=2$, which is verified by the simulation result. Actually, with $m=10^6$, the efficiency of $U_{oa}$ is $100.0\%$ when $t=3$ and $97.88\%$ when $t=2$.

\noindent{\large\bf{Examples for multi-sample and multi-dimensional cases.}}
Consider the multi-sample case. Suppose $d_1=d_2=2$, $n_1=n_2=9$ and the two samples are
\begin{eqnarray*}
X_6^{(1)}\leq X_8^{(1)}\leq X_2^{(1)}\leq X_4^{(1)}\leq X_7^{(1)}\leq X_5^{(1)}\leq X_3^{(1)}\leq X_9^{(1)}\leq X_1^{(1)}.\\
X_2^{(2)}\leq X_7^{(2)}\leq X_3^{(2)}\leq X_6^{(2)}\leq X_1^{(2)}\leq X_4^{(2)}\leq X_5^{(2)}\leq X_9^{(2)}\leq X_8^{(2)}.
\end{eqnarray*}
Then we have $L=3$ groups listed as
$G_1^{(1)}=\{6,8,2\},G_2^{(1)}=\{4,7,5\},G_3^{(1)}=\{3,9,1\}$ and $G_1^{(2)}=\{2,7,3\},G_2^{(2)}=\{6,1,4\},G_3^{(2)}=\{5,9,8\}$. An example of $OA(m=9,d=4,L=3,t=2)$ in step 1 is given as follows in transpose.
\begin{eqnarray*}
A^T=
\left(
\begin{array}{ccccccccc}
1&1&1&{ 2}&2&2&3&3&3\\
1&2&3&{ 1}&2&3&1&2&3\\
1&2&3&{ 2}&3&1&3&1&2\\
1&2&3&{ 3}&1&2&2&3&1
\end{array}
\right).
\end{eqnarray*}
Then we could possibly have the $\mathcal{X}_{\eta^i}$, $i=1,\ldots,9$, used in the construction of 9-run multi-sample construction as follows.
\begin{eqnarray*}
\{\mathcal{X}_{\eta^1},\ldots,\mathcal{X}_{\eta^9}\}=\left\{
\begin{array}{ccccccccc}
X^{(1)}_8&X^{(1)}_2&X^{(1)}_6&X^{(1)}_4&X^{(1)}_4&X^{(1)}_5&X^{(1)}_9&X^{(1)}_1&X^{(1)}_9\\
X^{(1)}_6&X^{(1)}_7&X^{(1)}_3&X^{(1)}_8&X^{(1)}_7&X^{(1)}_1&X^{(1)}_6&X^{(1)}_7&X^{(1)}_3\\
X^{(2)}_7&X^{(2)}_1&X^{(2)}_5&X^{(2)}_4&X^{(2)}_8&X^{(2)}_3&X^{(2)}_9&X^{(2)}_2&X^{(2)}_6\\
X^{(2)}_3&X^{(2)}_6&X^{(2)}_9&X^{(2)}_5&X^{(2)}_3&X^{(2)}_1&X^{(2)}_6&X^{(2)}_8&X^{(2)}_3
\end{array}
\right\}.
\end{eqnarray*}

Consider the multi-dimensional case. Suppose $X_1=(1.0,3.2)$, $X_2=(0.9,1.0)$, $X_3=(0.9,3.1)$, $X_4=(0.8,2.1)$, $X_5=(0.7,2.2)$, $X_6=(0.9,1.2)$, $X_7=(0.9,1.9)$, $X_8=(0.8,1.1)$, $X_9=(0.9,2.8)$. Simple clustering methods reveal $G_1=\{6,8,2\},G_2=\{4,7,5\},G_3=\{3,9,1\}$. The choosing of $\eta^i$, $i=1,\ldots,9$, might be the same as (\ref{eq:example}).

\end{document}